\documentclass[12pt]{article}
\usepackage{amsmath,amssymb}
\def\a{\alpha}
\def\b{\beta}

\def\8{\infty}
\def\sgn{\mbox{{\rm sgn}}}

\newtheorem {thm}{Theorem}[section]
\newtheorem {lem}[thm]{Lemma}
\newtheorem {cor}[thm]{Corollary}
\newtheorem {prop}[thm]{Proposition}

\makeatletter

\@addtoreset{equation}{section}
\makeatother
\begin{document}
\title{The $q$-Dixon--Anderson integral and \\ multi-dimensional $_1\psi_1$ summations
}
\author{
M{\small ASAHIKO} I{\small TO}%
\footnote{
School of Science and Technology for Future Life,
Tokyo Denki University,
Tokyo 120-8551, Japan 
\quad email: {\tt mito@cck.dendai.ac.jp}}
\ and P{\small ETER} J. F{\small ORRESTER}%
\footnote{
Department of Mathematics and Statistics, The University of Melbourne, Victoria 3010, Australia 
\quad email: {\tt p.forrester@ms.unimelb.edu.au}
}
}
\date{}
\maketitle
\begin{abstract}
\noindent
The Dixon--Anderson integral is a multi-dimensional integral evaluation fundamental to the theory of the Selberg integral.
The $_1\psi_1$ summation is a bilateral generalization of the $q$-binomial theorem. It is shown that a $q$-generalization of the
Dixon--Anderson integral, due to Evans, and multi-dimensional generalizations of the  $_1\psi_1$ summation, due to
Milne and Gustafson, can be viewed as having a common origin in the theory of $q$-difference equations as expounded by
Aomoto. Each is shown to be determined by a $q$-difference equation of rank one, and a certain asymptotic behavior.
In calculating the latter, essential use is made of the concepts of truncation, regularization and connection formulae.
\end{abstract}



\section{Introduction}
The Selberg integral is a multi-dimensional generalization of the Euler beta integral \cite{Se44}.
As surveyed in \cite{FW08}, after a dormant period of over thirty years since its discovery in the early 1940's, 
it shot into prominence upon the realization of its relevance to random matrix theory, the combinatorics of root systems,
and orthogonal polynomials in many variables, amongst other topics of contemporary interest.
In 1991 Anderson \cite{An91} gave a new proof of the Selberg integral by deriving a certain recurrence
in the number of integral variables $n$. At the heart of this proof was a further multi-dimensional integral
\begin{eqnarray}
\label{eq:DA}
&&\int_{z_n=x_{n-1}}^{x_{n}}\cdots\int_{z_2=x_1}^{x_2}\int_{z_1=x_0}^{x_1}
\prod_{i=1}^n\prod_{j=0}^{n}|z_i-x_j|^{s_j-1}\prod_{1\le k<l\le n}(z_l-z_k)\,dz_1dz_2\cdots dz_n
\nonumber\\
&&\quad =\frac{\Gamma(s_0)\Gamma(s_1)\cdots\Gamma(s_n)}{\Gamma(s_0+s_1+\cdots+s_n)}
\prod_{0\le i<j\le n}(x_j-x_i)^{s_i+s_j-1}.
\end{eqnarray}
For many years it was thought that (\ref{eq:DA}) was itself a new gamma function evaluation of a multi-dimensional
integral. However it was to transpire that (\ref{eq:DA}) in fact was first derived in a paper of Dixon \cite{Di05} 
written over eighty-five years earlier (see \cite{FW08} for the history of how Dixon's paper was rediscovered in the modern era).

One line of research generated by the Selberg integral was the study of $q$-generalizations. Askey \cite{As80}
conjectured the evaluation, in terms of $q$-gamma functions, of several multi-dimensional Jackson integrals which reduce
to the Selberg integral in the limit $q \to 1$. Evans  \cite{Ev92,Ev94} showed how two of these could be proved
by adopting the strategy of Anderson. This of course required a $q$-generalization of (\ref{eq:DA}). Evans \cite[Theorem 1, (2.5), p.759]{Ev92}
derived the sought $q$-generalization as
\begin{eqnarray}
\label{eq:Evans1}
&&\int_{z_n=x_{n-1}}^{x_{n}}\cdots\int_{z_2=x_1}^{x_2}\int_{z_1=x_0}^{x_1}
\prod_{i=1}^n\prod_{j=0}^{n}(qz_i/x_j;q)_{s_j-1}\prod_{1\le k<l\le n}(z_l-z_k)\,
\,d_qz_1d_qz_2\cdots d_qz_n
\nonumber\\
&&\quad =\frac{\Gamma_q(s_0)\Gamma_q(s_1)\cdots \Gamma_q(s_{n})}
{\Gamma_q(s_0+s_1+\cdots+s_{n})}
\prod_{0\le i<j\le n}x_j(x_i/x_j;q)_{s_j}(qx_j/x_i;q)_{s_i-1}.
\end{eqnarray}
We will refer to (\ref{eq:Evans1}) as the $q$-Dixon--Anderson integral.

Independent of the work of Anderson, Gustafson \cite{Gu90} invented the same strategy of using an auxiliary multi-dimensional integral
(nowadays referred to as a type I $q$-hypergeometric integral) 
to prove a product of $q$-gamma functions evaluation of a different generalization of 
the Selberg integral (type II $q$-hypergeometric integral relating to the $BC_n$ root system; 
the original Selberg integral relates to the $A_{n-1}$ root system). 
To implement this strategy Gustafson had to formulate and prove an appropriate $BC_n$ version of the
$q$-Dixon--Anderson integral \cite{Gu92}. But this was prior to the work of Evans, so (\ref{eq:Evans1}) was unknown. Rather the knowledge base of
Gustafson included his earlier work \cite{Gu87}, generalizing a result of Milne \cite{Mi85}, to give the product form evaluation of
the bilateral sum
\begin{equation}\label{MG}
\sum_{y_1,\dots,y_n = - \infty}^\infty
\prod_{1 \le i < j \le n} \Big ( \frac{z_i q^{y_i} - z_j q^{y_j} }{ z_i - z_j} \Big )
\prod_{i,j = 1}^n \frac{ (a_i z_j/z_i)_{y_j} }{(b_i z_j/z_i)_{y_j} }
t^{y_1 + \cdots + y_n}.
\end{equation}
As noted in \cite{Gu87}, the case $n=1$ of (\ref{MG}) is the definition of  the ${}_1 \psi_1$ series, and thus
the product form evaluation of (\ref{MG}) corresponds to a multi-dimensional generalization of the 
Ramanujan  ${}_1 \psi_1$ summation theorem. We will refer to (\ref{MG}) as the Milne--Gustafson summation.

In this paper we will show that (\ref{eq:Evans1}) and the product form evaluation of (\ref{MG}) are intimately related. The relationship
is seen by seeking an explanation for the product expressions from the viewpoint of $q$-difference systems.
One hint
of a common underpinning comes from (\ref{eq:Evans1}) and (\ref{MG}) both permitting generalizations involving Macdonald polynomials.
Thus it was pointed out in \cite{FR05} that the case $s_0 = \cdots = s_n = s$ 
of (\ref{eq:Evans1}) 
corresponds to  Okounkov's \cite{Ok98} integral
formula for the Macdonald polynomials $P_\kappa(x_0,\dots,x_n;q,q^s)$ in the case $\kappa = \emptyset$. On the other
hand one viewpoint of (\ref{MG}) is as a multi-dimensional ${}_1\psi_1$ summation associated to the root system $A_{n-1}$
(see e.g.~\cite{MS02} and references therein), and such generalized hypergeometric function identities  allow natural
generalizations to include Macdonald polynomials \cite{Kan96,BF99}.
\par
We remark that an elliptic analogue of the Milne--Gustafson summation is due to 
Kajihara--Noumi \cite{KN03}, and 
Rosengren \cite{Ro04,Ro06}, independently. In these references a specialization of (\ref{MG}) 
is generalized to involve elliptic analogues of the $q$-products, 
and the resulting summation is 
further more extended to an elliptic analogue of Kajihara's $q$-transformation identity \cite{Kaj04} 
between a summation over $n$ variables and a summation over $m$ variables.
This gives another hint of a relationship between (\ref{eq:Evans1}) and the product form evaluation of (\ref{MG}). 
Thus in the original paper of Dixon \cite{Di05} the Dixon--Anderson integral (\ref{eq:DA}) 
is generalized to a transformation identity between an $n$ and an $m$ dimensional multi-dimensional integral. 
And in proving a conjectured elliptic generalization of the
Selberg integral due to van Diejen and Spiridonov \cite{DS01}, Rains \cite{Ra10} has proved a transformation identity between
multi-dimensional elliptic integrals, which he has subsequently \cite{Ra09} shown permits under certain limiting operations
a reduction to the original transformation identity of Dixon.
\par
We will discuss three multi-dimensional bilateral series, 
\begin{enumerate}
\item[(1)] the Milne--Gustafson summation formula (Theorem \ref{thm:MG}),
\item[(2)] a multi-dimensional bilateral extension of Evans's $q$-Dixon--Anderson integral (Theorem \ref{thm:02}), 
\item[(3)] the product expression for Gustafson's $A_n$ sum contained in \cite{Gu87}, itself generalizing
the Milne--Gustafson summation (Theorem \ref{thm:03main}).
\end{enumerate}
\noindent
The sum of the case (1) connects the sums of (2) and (3) 
as a hub. As already mentioned, the aim of this paper is to give an explanation for the formulae of product expression of these sums 
from a view-point of $q$-difference equations. 

The method for proving these results is consistent with the concepts introduced by Aomoto 
and Aomoto--Kato in the early 1990's
in the series of papers \cite{Ao90,Ao91,Ao94,Ao95-1,AK91,AK93,AK94-1,AK94-2}. 
Aomoto showed an isomorphism between a class of the Jackson integrals of hypergeometric type,
which he called the {\it $q$-analog de Rham cohomology} \cite{Ao90,Ao91}, and a class of theta functions, i.e., 
holomorphic functions possessing a quasi-periodicity \cite[Theorem 1]{Ao95-1}. 
This isomorphism indicates that it is essential to analyze 
both in order to know the structure of $q$-hypergeometric functions, in particular, 
the meaning of known special formulae. 
In this paper the process to obtain the holomorphic functions through this isomorphism is called {\it regularization}. 
When we fix a basis of the class of holomorphic functions as a linear space, 
an arbitrary function of the space can be expressed as a linear combination of the elements of the specific basis, 
which he called the {\it connection formula} \cite[Theorem 3]{Ao95-1}. 
(As its simplest examples, Ramanujan's $_1\psi_1$ summation formula and 
the $q$-Selberg integral \cite{As80,Ha88,Kad88,Ev92}  have been explained. See \cite[Examples 1, 2]{Ao95-1}.) 
One way to choose a good basis is through its asymptotic behavior upon a limiting process with  
respect to parameters included in the definition of the Jackson integral of hypergeometric type. %
Moreover the asymptotic behavior can be calculated from Jackson integrals possessing 
appropriate cycles which include their critical points. We call the process to fix the cycles the {\it truncation}.
(These cycles are called the {\it characteristic cycles} \cite{AK94-2} 
or the {\it $\alpha$-stable or $\alpha$-unstable cycles} \cite{Ao94} by Aomoto. 
The meaning of ``$\alpha$" is mentioned in Section \ref{section:01}. 
The word {\it truncation} itself is first used by van Diejen in another context \cite{vD97, Ito06-2}.)
There is another way to characterize the connection formula.  
This is by showing that 
a multi-dimensional bilateral series originating as a general solution 
of the $q$-difference equation of the Jackson integrals with respect to parameters 
is expressed as a linear combination of multi-dimensional unilateral series
as special solutions, each fixed by their asymptotic behaviors \cite[Theorem (4.2)]{Ao94}. 
(We can see different examples of $q$-difference equations and the connection formulae from this paper 
in \cite{Ito08,Ito09,IS08}, and \cite{IS08} explains 
the Sears--Slater transformation for the very-well-poised $q$-hypergeometric series 
from the view-point of this paper in the setting of $BC$ type symmetry. 
In the recent survey article \cite{IF13} we have detailed this program as it
applies to Ramanujan's ${}_1 \psi_1$ summation formula and Bailey's very-well-poised ${}_6 \psi_6$ summation,
and show how these logically lead to the consideration of higher dimensional extensions such as (2)). 
We discuss three product formulae for the Jackson integrals corresponding to the sums of (1), (2) and (3)
as simple examples of this concept.

This paper is organized as follows. 
After defining basic terminology in Section \ref{section:00}, 
we first show the product expression of the Milne--Gustafson sum 
using concepts of {\it truncation}, {\it regularization} and {\it connection formulae} in Section \ref{section:01}. 
Though the Milne--Gustafson sum can be obtained from our other two examples, 
we explain it individually, because the Milne--Gustafson sum has simpler structure than the other two,
and it is instructive in outlining the concepts of this paper. 
The subsequent two sections are devoted to explaining Evans's $q$-Dixon--Anderson integral   
and Gustafson's $A_n$ sum, respectively. 
Although technically more involved than the Milne--Gustafson sum, the overall strategy will
 be seen to be the same. 
Specifically in section \ref{section:02} we introduce a bilateral extension of Evans's $q$-Dixon--Anderson integral 
and use the method of $q$-difference equations to deduce its evaluation. 
In section \ref{section:03}, for Gustafson's $A_n$ sum we also use the method of $q$-difference equations to deduce its evaluation. 
In reading sections 3 to 5, a repetition of the main steps needed to implement the $q$-difference
equation method will become apparent. \\

\section{Definition of the Jackson integral}
\label{section:00}
Throughout this paper, we fix $q$ as $0<q<1$ and use the symbols 
$(a)_\8:=\prod_{i=0}^\8(1-q^i a)$ and $(a)_N:=(a)_\8/(q^Na)_\8$. 
We define $\theta(a)$ by $\theta(a):=(a)_\8(q/a)_\8$, 
which satisfies 
\begin{equation}
\label{eq:00quasi-period}
\theta(qa)=-\theta(a)/a.
\end{equation}
Let $S_n$ be the symmetric group on $\{1,2,\ldots, n\}$.
For a function $f(z)=f(z_1,z_2,\ldots,z_n)$ on $(\mathbb{C}^*)^n$, 
we define action of the symmetric group $S_n$ on $f(z)$ by 
$$
(\sigma f)(z):=f(\sigma^{-1}(z))=f(z_{\sigma(1)},z_{\sigma(2)},\ldots,z_{\sigma(n)})
\quad\mbox{for}\quad \sigma\in S_n.
$$
We say that a function $f(z)$ on $(\mathbb{C}^*)^n$
is {\it symmetric} or {\it skew-symmetric} 
if $\sigma f(z)=f(z)$ or $\sigma f(z)=(\sgn\,\sigma )\,f(z)$ 
for all $\sigma \in S_n$, respectively.
We denote by ${\cal A} f(z)$ 
the alternating sum over $S_n$ defined by 
\begin{equation}
\label{eq:00Af}
({\cal A} f)(z):=\sum_{\sigma\in S_n}(\sgn\, \sigma)\,\sigma f(z),
\end{equation}
which is skew-symmetric. 
\par
For $a,b\in \mathbb{C}$, we define 
\begin{equation}
\label{eq:00jac1}
\int_a^b f(z)d_qz:=\int_0^b f(z)d_qz-\int_0^a f(z)d_qz,
\end{equation}
where
$$
\int_0^a f(z)d_qz:=(1-q)\sum_{\nu=0}^\8 f(aq^\nu)aq^\nu,
$$
which is called the {\it Jackson integral}. As $q\to 1$, 
$\int_a^b f(z)d_qz\to \int_a^b f(z)dz$ \cite{AAR99}. 
In this paper we basically use the Jackson integral of multiplicative measure as 
$$
\int_0^a f(z)\frac{d_qz}{z}=(1-q)\sum_{\nu=0}^\8 f(aq^\nu). 
$$
Let $\mathbb{N}$ be the set of non-negative integers. 
For a function $f(z)=f(z_1,\ldots,z_n)$ on $(\mathbb{C}^*)^n$ and 
an arbitrary point $x=(x_1,\ldots,x_n)\in (\mathbb{C}^*)^n$, 
we define the multiple Jackson integral as
\begin{equation}
\label{eq:00jac2}
\int_0^{\mbox{\small $x$}}f(z)\,\frac{d_qz_1}{z_1}\wedge\cdots\wedge\frac{d_qz_n}{z_n}
:=(1-q)^n\sum_{(\nu_1,\ldots,\nu_n)\in {\mathbb{N}}^n}f(x_1 q^{\nu_1},\ldots,x_n q^{\nu_n}). 
\end{equation}
In this paper we use the multiple bilateral sum extending the above Jackson integral
\begin{equation}
\label{eq:00jac3}
\int_0^{\mbox{\small $x$}\8}f(z)\,\frac{d_qz_1}{z_1}\wedge\cdots\wedge\frac{d_qz_n}{z_n}
:=(1-q)^n\sum_{(\nu_1,\ldots,\nu_n)\in {\mathbb{Z}}^n}f(x_1 q^{\nu_1},\ldots,x_n q^{\nu_n}),
\end{equation}
which we also call the {\it Jackson integral}. By definition the Jackson integral (\ref{eq:00jac3}) 
is invariant under the shift $x_i\to qx_i$, $1\le i\le n$.
While we can consider the limit $q\to1$ for the Jackson integral (\ref{eq:00jac2}) defined over $\mathbb{N}^n$,  
the Jackson integral (\ref{eq:00jac3}) defined over $\mathbb{Z}^n$ generally diverges if $q\to1$. 
However, as we will see later, since the {\it truncation} of 
the Jackson integral (\ref{eq:00jac3}) is corresponding to the sum (\ref{eq:00jac2}) over $\mathbb{N}^n$,  
if we need to consider the limit $q\to1$, we switch from (\ref{eq:00jac3}) to (\ref{eq:00jac2}) by the process of the truncation. 

We will state one of the key lemmas of this paper for deriving $q$-difference equations. 
For this let $\Phi(z)$ be a symmetric function on $(\mathbb{C}^*)^n$
and  for a function $\varphi(z)$, define the function $\nabla_{\!i}\varphi(z)$ ($1\le i\le n$) by
\begin{equation}
\label{eq:00nabla}
(\nabla_{\!i}\varphi)(z):=\varphi(z)-\frac{T_{z_i}\Phi(z)}{\Phi(z)}T_{z_i}\varphi(z),
\end{equation}
where $T_{z_i}$ means the shift operator of $z_i\to qz_i$, i.e.,
$T_{z_i}f(\ldots,z_i,\ldots)=f(\ldots,qz_i,\ldots)$. 
We then have  
\begin{lem}
\label{lem:00nabla=0}
For a meromorphic function $\varphi(z)$ on $(\mathbb{C}^*)^n$, if the integral 
$$
\int_0^{\mbox{\small $x$}\8}\varphi(z)\Phi(z)\frac{d_qz_1}{z_1}\wedge\cdots\wedge\frac{d_qz_n}{z_n}
$$
converges, then 
\begin{equation}
\label{eq:00nabla=0}
\int_0^{\mbox{\small $x$}\8}\Phi(z)\nabla_{\!i}\varphi(z)\frac{d_qz_1}{z_1}\wedge\cdots\wedge\frac{d_qz_n}{z_n}=0. 
\end{equation}
Moreover, 
\begin{equation}
\label{eq:00A}
\int_0^{\mbox{\small $x$}\8}\Phi(z){\cal A}\nabla_{\!i}\varphi(z)\frac{d_qz_1}{z_1}\wedge\cdots\wedge\frac{d_qz_n}{z_n}=0,
\end{equation}
where ${\cal A}$ indicates the skew-symmetrization defined in {\rm (\ref{eq:00Af})}. 
\end{lem}
{\bf Proof.} 
From the definition (\ref{eq:00nabla}) of $\nabla_{\!i}$, 
(\ref{eq:00nabla=0}) is equivalent to the statement 
$$
\int_0^{\mbox{\small $x$}\8}\varphi(z)\Phi(z)\frac{d_qz_1}{z_1}\wedge\cdots\wedge\frac{d_qz_n}{z_n}
=\int_0^{\mbox{\small $x$}\8}T_{z_i}\varphi(z)\, T_{z_i}\Phi(z)
\frac{d_qz_1}{z_1}\wedge\cdots\wedge\frac{d_qz_n}{z_n},
$$
if the left-hand side converges. Moreover this equation 
is just confirmed from the fact that the Jackson integral is invariant under the $q$-shift $z_i\to qz_i$ ($1\le i\le n$).
Next we will confirm (\ref{eq:00A}). 
Using 
$
\sigma\Phi(z)=\Phi(z)
$,
we have 
\begin{eqnarray*}
\Phi(z){\cal A}\nabla_{\!i}\varphi(z)
&=&\Phi(z)\sum_{\sigma\in S_n}(\sgn\, \sigma)\,\sigma(\nabla_{\!i}\varphi)(z)
=\sum_{\sigma\in S_n}(\sgn\, \sigma)\,\sigma\Phi(z)\sigma(\nabla_{\!i}\varphi)(z)\\
&=&\sum_{\sigma\in S_n}(\sgn\, \sigma)\,\sigma\Big(\Phi(z)\nabla_{\!i}\varphi(z)\Big),
\end{eqnarray*}
so that we obtain 
$$
\int_0^{\mbox{\small $x$}\8}\Phi(z){\cal A}\nabla_{\!i}\varphi(z)\frac{d_qz_1}{z_1}\wedge\cdots\wedge\frac{d_qz_n}{z_n}
=\sum_{\sigma\in S}(\sgn\, \sigma)
\int_0^{\mbox{\small $x$}\8}\sigma\Big(\Phi(z)\nabla_{\!i}\varphi(z)\Big)\frac{d_qz_1}{z_1}\wedge\cdots\wedge\frac{d_qz_n}{z_n}
$$
$$
=\sum_{\sigma\in S}(\sgn\, \sigma)
\int_0^{\sigma^{-1}\mbox{\small $x$}\8}\!\!\!\!\!
\Phi(z)\nabla_{\!i}\varphi(z)\frac{d_qz_1}{z_1}\wedge\cdots\wedge\frac{d_qz_n}{z_n}
=\sum_{\sigma\in S}(\sgn\, \sigma)\,\sigma\!\!
\int_0^{\mbox{\small $x$}\8}\Phi(z)\nabla_{\!i}\varphi(z)\frac{d_qz_1}{z_1}\wedge\cdots\wedge\frac{d_qz_n}{z_n},
$$
which vanishes from (\ref{eq:00nabla=0}). $\square$
\section{Jackson integral of Milne--Gustafson type}
\label{section:01}
In this section, we will explain the methods of this paper using the Milne--Gustafson summation formula.
\subsection{Definitions and the result}
Let $a_1,\ldots,a_n$, $b_1,\ldots,b_n$ and $\alpha$ be complex numbers satisfying 
\begin{equation}
\label{eq:01condition01}
|qa_1^{-1}\cdots a_n^{-1}b_1^{-1}\cdots b_n^{-1}|<|q^\alpha|<1.
\end{equation}
Let $\Phi(z)$ and $\Delta(z)$ be the functions defined by 
\begin{equation}
\label{eq:01Phi1}
\Phi(z):=(z_1z_2\cdots z_n)^\alpha\prod_{i=1}^n\prod_{j=1}^{n}\frac{(qa_j^{-1}z_i)_\8}{(b_jz_i)_\8}
\end{equation}
and
\begin{equation}
\label{eq:01Delta1}
\Delta(z):=\prod_{1\le i<j\le n}(z_j-z_i).
\end{equation}
For $x=(x_1,x_2,\ldots,x_n)\in (\mathbb{C}^*)^n$, we define the sum $I(x)$ by  
\begin{equation}
\label{eq:01I(x)1}
I(x):=\int_0^{\mbox{\small $x$}\8}\Phi(z)\Delta(z)\,\frac{d_qz_1}{z_1}\wedge\cdots\wedge\frac{d_qz_n}{z_n},
\end{equation}
which converges 
absolutely 
under the condition (\ref{eq:01condition01}) 
(see \cite[Lemma 3.19]{Gu87}, \cite[Lemma 2.5]{Mi86}). 
We call $I(x)$ the {\it Jackson integral of Milne--Gustafson type}. 
By definition $I(x)$ is skew-symmetric. 
\begin{thm}[Milne--Gustafson \cite{Gu87,Mi85}]
\label{thm:MG}
 For an arbitrary $x\in (\mathbb{C}^*)^n$, $I(x)$ is evaluated as 
\begin{eqnarray}
\label{eq:01I(x)}
I(x)&=&
(1-q)^n
\frac{(q)_\8^n 
\prod_{i=1}^{n}\prod_{j=1}^{n}(qa_i^{-1}b_j^{-1})_\8
}
{(q^\alpha)_\8(q^{1-\alpha} a_1^{-1}a_2^{-1}\cdots a_{n}^{-1}b_1^{-1}b_2^{-1}\cdots b_{n}^{-1})_\8
}\nonumber\\
&&\times (x_1x_2\cdots x_n)^{\alpha}
\frac{\theta(q^{\alpha}x_1x_2\cdots x_n b_1b_2\cdots b_n)
}{\prod_{i=1}^n\prod_{j=1}^n\theta(x_ib_j)}\prod_{1\le i<j\le n}x_j\theta(x_i/x_j).
\end{eqnarray}
\end{thm}
{\bf Remark.} If $n=1$, (\ref{eq:01I(x)}) coincides with 
Ramanujan's $_1\psi_1$ summation theorem. In this sense, 
(\ref{eq:01I(x)}) is a natural multi-dimensional extension of the $_1\psi_1$ sum. 
On the other hand, 
if $q\to 1$, the limiting formula of a special case of (\ref{eq:01I(x)}) 
coincides with the Dixon--Anderson integral (\ref{eq:DA}) with $x_0=0$, 
as we will see in Remark after Corollary \ref{cor:01I(a)}. In this sense, 
the Milne--Gustafson sum is a natural extension of both quantities. 
\\

The rest of this section is devoted to explaining the Milne--Gustafson summation formula (\ref{eq:01I(x)}) 
as a connection formula between solutions of 
a $q$-difference equation, and consequently this leads us to a simple proof.  
\subsection{$q$-difference equation with respect to $\alpha$}
In this subsection we derive the $q$-difference equation with respect to $\alpha$, which $I(x)$ satisfies. 
We use $I(\alpha;x)$ instead of $I(x)$ to see the $\alpha$ dependence. 
\begin{lem}
\label{lem:rec}
The recurrence relation for $I(\alpha;x)$ is given by 
\begin{equation}
\label{eq:01rec1}
I(\alpha;x)=\frac{1-q^{\alpha}a_1a_2\cdots a_nb_1b_2\cdots b_n}{a_1a_2\cdots a_n(1-q^\alpha)}I(\alpha+1;x).
\end{equation}
\end{lem}
{\bf Proof.}
Since the ratio ${T_{z_1}\Phi(z)}/{\Phi(z)}$ is written as 
$$
\frac{T_{z_1}\Phi(z)}{\Phi(z)}=q^\alpha\prod_{j=1}^n\frac{1-b_jz_1}{1-q a_j^{-1}z_1},
$$
if we put  
$
\varphi(z)=z_2^{n-1}z_3^{n-2}\cdots z_n\prod_{j=1}^n(1- a_j^{-1}z_1),
$
then, from (\ref{eq:00nabla}) we have 
$$
\nabla_{\!1}\varphi(z)=z_2^{n-1}z_3^{n-2}\cdots z_n\Big(
\prod_{j=1}^n(1- a_j^{-1}z_1)
-q^\alpha\prod_{j=1}^n(1- b_jz_1)\Big),
$$
so that the skew-symmetrization of the above equation is given by 
\begin{equation}
\label{eq:01Anabla}
{\cal A}\nabla_{\!1}\varphi(z)=\Big((-1)^{n-1}(1-q^\alpha)+(-1)^n
\frac{1-q^{\alpha}a_1a_2\cdots a_nb_1b_2\cdots b_n}{a_1a_2\cdots a_n}z_1z_2\cdots z_n\Big)\Delta(z).
\end{equation}
Since 
$$
I(\alpha+1;x)=
\int_0^{\mbox{\small $x$}\8}z_1z_2\cdots z_n
\Phi(z)\Delta(z)\,\frac{d_qz_1}{z_1}\wedge\cdots\wedge\frac{d_qz_n}{z_n},
$$
using (\ref{eq:00A}) of Lemma \ref{lem:00nabla=0} for (\ref{eq:01Anabla}), we obtain the relation (\ref{eq:01rec1}). $\square$

\subsection{Truncation}
For the special point $x=a=(a_1,a_2,\ldots,a_n)$, we call $I(a)$ the 
{\it truncated Jackson integral of Milne--Gustafson type}. By definition 
$I(a)$ is a sum over $\mathbb{N}^n$, while $I(x)$ is generally a sum over the lattice $\mathbb{Z}^n$.
It has the advantage of simplifying the computation of the $\alpha\to +\8$ asymptotic behavior,
as will be seen below.
(The lattice  $\{(x_1q^{\nu_1},\ldots,x_nq^{\nu_n})\in (\mathbb{C}^*)^n\,;\,(\nu_1,\ldots,\nu_n)\in \mathbb{Z}^n\}$ is 
called the {\it $q$-cycle} \cite{Ao91} of $I(x)$, while the set $\{(a_1q^{\nu_1},\ldots,a_nq^{\nu_n})\in (\mathbb{C}^*)^n\,;\,(\nu_1,\ldots,\nu_n)\in \mathbb{N}^n\}$ as the support of $I(a)$ is called 
the {\it $\alpha$-stable cycle} in \cite{Ao94,AK94-2}.)
\begin{lem}
\label{lem:I(a+N;a)}
The asymptotic behavior of $I(\alpha+N;a)$ as $N\to +\8$ is given by
\begin{equation}
\label{eq:01asym1}
I(\alpha+N;a)\sim (1-q)^n(a_1a_2\cdots a_n)^{\alpha+N}\Delta(a)
\prod_{i=1}^n\prod_{j=1}^n\frac{(qa_i/a_j)_\8}{(a_ib_j)_\8}
\quad(N\to +\8).
\end{equation}
\end{lem}
{\bf Proof.} Since $(q^{1+\nu_i})_\8=0$ if $\nu_i<0$, by definition $I(\alpha+N;a)$ is written as 
\begin{eqnarray}
\label{eq:01asym1.5}
I(\alpha+N;a)&=&(1-q)^n\!\!\!\!\!
\sum_{(\nu_1,\ldots,\nu_n)\in \mathbb{N}^n}
(a_1a_2\cdots a_nq^{\nu_1+\cdots+\nu_n})^{\alpha+N}\nonumber\\
&&\times\prod_{i=1}^n\prod_{j=1}^{n}\frac{(a_j^{-1}a_iq^{1+\nu_i})_\8}{(b_ja_iq^{\nu_i})_\8}
\prod_{1\le i<j\le n}(a_jq^{\nu_j}-a_iq^{\nu_i}).
\end{eqnarray}
Since $(p)_\8\to 1$ if $p\to 0$, 
the term $|\prod_{i=1}^n\prod_{j=1}^{n}(a_j^{-1}a_iq^{1+\nu_i})_\8/(b_ja_iq^{\nu_i})_\8|$ is bounded for 
$(\nu_1,\ldots,\nu_n)\in \mathbb{N}^n$ and $N>0$. 
Then the leading term of the asymptotic behavior of $I(\alpha+N;a)$ as $N\to +\8$ 
depends on the maximum value of 
$|(a_1a_2\cdots a_nq^{\nu_1+\cdots+\nu_n})^{\alpha+N}
\prod_{1\le i<j\le n}(a_jq^{\nu_j}-a_iq^{\nu_i})|$
over $(\nu_1,\ldots,\nu_n)\in \mathbb{N}^n$. 
Since $|\prod_{i=1}^n\prod_{j=1}^{n}(a_j^{-1}a_iq)_\8/(b_ja_i)_\8|\ne 0$, 
the summand in the right-hand side of (\ref{eq:01asym1.5}) corresponding to 
$(\nu_1,\ldots,\nu_n)=(0,\ldots,0)$ gives the leading term (\ref{eq:01asym1}).
$\square$ \\[2pt]
\indent
From this lemma, with repeated use of Lemma \ref{lem:rec}, we immediately have  
\begin{cor}
\label{cor:01I(a)}
The truncated Jackson integral $I(a)$ is expressed as  
\begin{equation}
\label{eq:01I(a)}
I(a)=
(1-q)^n(a_1a_2\cdots a_n)^{\alpha}
\frac{(q)_\8^n 
(q^\alpha a_1a_2\cdots a_{n}b_1b_2\cdots b_{n})_\8
}
{(q^\alpha)_\8\prod_{i=1}^{n}\prod_{j=1}^{n}(a_ib_j)_\8}
\prod_{1\le i<j\le n}a_j\theta(a_i/a_j).
\end{equation}
\end{cor}
{\bf Proof.} 
By repeated use of the recurrence relation (\ref{eq:01rec1}), we have 
$$I(\alpha;x)
=\frac{(q^{\alpha}a_1a_2\cdots a_nb_1b_2\cdots b_n)_N}{(a_1a_2\cdots a_n)^N(q^\alpha)_N}I(\alpha+N;x).$$
If we put $x=a$ and take $N\to +\8$, we obtain  
$$
I(\alpha;a)
=\frac{(q^{\alpha}a_1a_2\cdots a_nb_1b_2\cdots b_n)_\8}{(q^\alpha)_\8}\lim_{N\to \8}\frac{I(\alpha+N;a)}{(a_1a_2\cdots a_n)^N}, 
$$
which coincides with the right-hand side of (\ref{eq:01I(a)}) if we use (\ref{eq:01asym1}).
$\square$\\

\noindent
{\bf Remark.} If $n=1$, (\ref{eq:01I(a)}) coincides with the $q$-binomial theorem which is a special case of 
Ramanujan's $_1\psi_1$ summation theorem. 
The truncated Jackson integral $I(a)$ is the sum over $\mathbb{N}^n$ by definition, 
and the definition of this case coincides with that of the ordinary Jackson integral explained by (\ref{eq:00jac2}). 
Therefore we can consider the limiting case of (\ref{eq:01I(a)}) as $q\to 1$. 
If we substitute $\alpha$, $a_j$ and $b_j$ as 
$\alpha\to s_0$, $a_j\to x_{j-1}$ and $b_j\to q^{s_{j-1}}/x_{j-1}$, respectively, 
then (\ref{eq:01I(a)}) in the limit $q\to 1$ coincides with the Dixon--Anderson integral (\ref{eq:DA}) with $x_0=0$. 

\subsection{Regularization and connection formulae}
Let ${\cal I}(x)$ and $h(x)$ be the functions defined by 
\begin{equation}
\label{eq:01cal I(x)1}
{\cal I}(x)=\frac{I(x)}{h(x)},
\quad 
h(x)=
(x_1x_2\cdots x_n)^{\alpha}
\frac{\prod_{1\le i<j\le n}x_j\theta(x_i/x_j)}
{\prod_{i=1}^{n}\prod_{j=1}^{n}\theta(b_jx_i)}.
\end{equation}
We call ${\cal I}(x)$ the {\it regularized Jackson integral} of $I(x)$. 
Since the trivial poles and zeros of $I(x)$ are canceled out by multiplying together $1/h(x)$ and $I(x)$, 
we have the following. 
\begin{lem}
\label{lem:h-and-s}
The regularization ${\cal I}(x)$ is holomorphic on $(\mathbb{C}^*)^n$ and symmetric. 
\end{lem}
\noindent
{\bf Proof.} From the expression (\ref{eq:01Phi1}) of $\Phi(z)$ as integrand of (\ref{eq:01I(x)1}), 
the function $I(x)$ has the poles lying only in the set 
$\{x=(x_1,x_2,\ldots,x_n)\in (\mathbb{C}^*)^n\,;\,\prod_{i=1}^n\prod_{j=1}^n\theta(b_jx_i)=0
\}$.
Moreover $I(x)$ is divisible by $x_j\theta(x_i/x_j)$ 
because $I(x)$ is skew-symmetric 
and invariant under $q$-shift
. We therefore obtain 
$$
I(x)={\cal I}(x)h(x),
$$
where ${\cal I}(x)$ is some holomorphic function on $(\mathbb{C}^*)^n$. 
Since $I(x)$ and $h(x)$ are both skew-symmetric, ${\cal I}(x)$ is symmetric. 
$\square$\\

From (\ref{eq:01cal I(x)1}) and the quasi-periodicity (\ref{eq:00quasi-period}) of the theta function, we have 
$$\frac{T_{x_i}h(x)}{h(x)}=-q^\alpha x_1x_2\cdots x_nb_1b_2\cdots b_n.$$
And since $I(x)$ is an invariant under the $q$-shift $x_i\to qx_i$, the holomorphic function ${\cal I}(x)$ on $(\mathbb{C}^*)^n$ 
satisfies the $q$-difference equations 
\begin{equation}
\label{eq:01quasi-period}
T_{x_i}{\cal I}(x)=-\frac{{\cal I}(x)}{q^\alpha x_1x_2\cdots x_nb_1b_2\cdots b_n},\quad i=1,2,\ldots,n.
\end{equation}
Since the set of holomorphic functions on $(\mathbb{C}^*)^n$ satisfying (\ref{eq:01quasi-period}) 
has the dimension 1 as a $\mathbb {C}$-linear space, 
we can take $\theta(q^\alpha x_1x_2\cdots x_nb_1b_2\cdots b_n)$ as a basis of the linear space. Thus 
${\cal I}(x)$ is uniquely expressed as 
\begin{equation}
\label{eq:01cal I(x)2}
{\cal I}(x)=C\,\theta(q^\alpha x_1x_2\cdots x_nb_1b_2\cdots b_n),
\end{equation}
where $C$ is some constant independent of $x$. 
\begin{lem}[connection formula]
For arbitrary $x,y\in (\mathbb{C}^*)^n$, the {\it connection formula} between ${\cal I}(x)$ and ${\cal I}(y)$ 
is written as 
\begin{equation}
\label{eq:01cal I(x)I(y)}
{\cal I}(x)=\frac{\theta(q^\alpha x_1x_2\cdots x_nb_1b_2\cdots b_n)}
{\theta(q^\alpha y_1y_2\cdots y_nb_1b_2\cdots b_n)}
{\cal I}(y).
\end{equation}
In particular, if we set $y=a\in (\mathbb{C}^*)^n$, then 
\begin{equation}
\label{eq:01cal I(x)I(a)}
{\cal I}(x)=\frac{\theta(q^\alpha x_1x_2\cdots x_nb_1b_2\cdots b_n)}
{\theta(q^\alpha a_1a_2\cdots a_nb_1b_2\cdots b_n)}
{\cal I}(a).
\end{equation}
\end{lem}
{\bf Proof.}
From (\ref{eq:01cal I(x)2}), we immediately have (\ref{eq:01cal I(x)I(y)}). $\square$\\[6pt]
{\bf Remark.} 
If we switch the symbols from ${\cal I}(x)$ to $I(x)$, we have 
$$
{I}(x)=\frac{h(x)\theta(q^\alpha x_1x_2\cdots x_nb_1b_2\cdots b_n)}
{h(y)\theta(q^\alpha y_1y_2\cdots y_nb_1b_2\cdots b_n)}
{I}(y).
$$
In particular, if we set $y=a$ in the above equation, we obtain 
\begin{equation}
\label{eq:01I(x)I(a)}
{I}(x)=\frac{h(x)\theta(q^\alpha x_1x_2\cdots x_nb_1b_2\cdots b_n)}
{h(a)\theta(q^\alpha a_1a_2\cdots a_nb_1b_2\cdots b_n)}
{I}(a),
\end{equation}
which is also the connection formula between a solution $I(x)$ of the $q$-difference equation (\ref{eq:01rec1}) 
and the special solution $I(a)$ fixed by its 
asymptotic behavior (\ref{eq:01asym1}) as $\alpha\to +\8$. In addition, its connection coefficient is written as 
a ratio of theta functions (i.e., that of $q$-gamma functions), and is of course invariant under the shift 
$\alpha\to \alpha+1$.  
Using the evaluation (\ref{eq:01I(a)}) of $I(a)$, 
the connection formula (\ref{eq:01I(x)I(a)}) exactly coincides with (\ref{eq:01I(x)}) for the Milne--Gustafson sum.\\

Using (\ref{eq:01I(a)}), the constant $C$ in (\ref{eq:01cal I(x)2}) is also calculated explicitly as 
\begin{equation}
\label{eq:01C}
C=\frac{{I}(a)}
{h(a)\theta(q^\alpha a_1a_2\cdots a_nb_1b_2\cdots b_n)}
=
\frac{(1-q)^n(q)_\8^n 
\prod_{i=1}^{n}\prod_{j=1}^{n}(qa_i^{-1}b_j^{-1})_\8
}
{(q^\alpha)_\8(q^{1-\alpha} a_1^{-1}a_2^{-1}\cdots a_{n}^{-1}b_1^{-1}b_2^{-1}\cdots b_{n}^{-1})_\8
}.
\end{equation}
From (\ref{eq:01cal I(x)2}) we therefore obtain 
\begin{prop}
For an arbitrary $x\in (\mathbb{C}^*)^n$, ${\cal I}(x)$ is expressed as
\begin{equation}
\label{eq:01cal I(x)3}
{\cal I}(x)=
(1-q)^n
\frac{(q)_\8^n 
\prod_{i=1}^{n}\prod_{j=1}^{n}(qa_i^{-1}b_j^{-1})_\8
}
{(q^\alpha)_\8(q^{1-\alpha} a_1^{-1}\cdots a_{n}^{-1}b_1^{-1}\cdots b_{n}^{-1})_\8
}\, 
\theta(q^{\alpha}x_1\cdots x_n b_1\cdots b_n).
\end{equation}
\end{prop}
Since (\ref{eq:01I(x)}) is equivalent to (\ref{eq:01cal I(x)3}), the Milne--Gustafson summation formula 
as the expression (\ref{eq:01cal I(x)3}) gives the isomorphism between 
the sets of $I(x)$ and of theta functions defined by (\ref{eq:01quasi-period}). \\

If we set 
\begin{equation}
\label{eq:01beta}
\beta:=1-\alpha_1-\cdots-\alpha_n-\beta_1-\cdots-\beta_n-\alpha, 
\end{equation}
where $\alpha_i$ and $\beta_i$ are given by $a_i=q^{\alpha_i}, b_i=q^{\beta_i}$, 
after rearrangement, 
the formula (\ref{eq:01cal I(x)3}) is also expressed as the following Macdonald-type sum, 
whose value is given by an $x$-independent constant 
\cite{Ma03,vD97,Ito06-1}.
\begin{prop} Under the condition $a_1\cdots a_n b_1\cdots b_nq^{\alpha+\beta}=q$, 
\begin{eqnarray}
\label{eq:01cal I(x)4}
&&\int_0^{\mbox{\small $x$}\8}\frac{\prod_{i=1}^n\prod_{j=1}^n (qa_j^{-1}z_i)_\8(qb_j^{-1}z_i^{-1})_\8}
{(q^{\beta}a_1\cdots a_n z_1^{-1}\cdots z_n^{-1})_\8(q^{\alpha}b_1\cdots b_n z_1\cdots z_n)_\8
\prod_{1\le i<j\le n}(qz_i/z_j)_\8(qz_j/z_i)_\8}\nonumber\\
&&
\qquad\qquad
\times\frac{d_qz_1}{z_1}\wedge\cdots\wedge\frac{d_qz_n}{z_n}
=\frac{(1-q)^n(q)_\8^n 
\prod_{i=1}^{n}\prod_{j=1}^{n}(qa_i^{-1}b_j^{-1})_\8
}
{(q^\alpha)_\8(q^\beta)_\8
}.
\end{eqnarray}
\end{prop} 
{\bf Proof.} Since $h(x)\theta(q^{\alpha}x_1\cdots x_n b_1\cdots b_n)$ is invariant under the $q$-shift $x_i\to qx_i$, 
from (\ref{eq:01cal I(x)2}), we have 
\begin{equation*}
\int_0^{\mbox{\small $x$}\8}\frac{\Phi(z)\Delta(z)}
{h(z)\theta(q^{\alpha}z_1\cdots z_n b_1\cdots b_n)}\,\frac{d_qz_1}{z_1}\wedge\cdots\wedge\frac{d_qz_n}{z_n}
=C,
\end{equation*}
so that 
\begin{equation*}
\int_0^{\mbox{\small $x$}\8}\frac{\prod_{i=1}^n\prod_{j=1}^n (qa_j^{-1}z_i)_\8(qb_j^{-1}z_i^{-1})_\8}
{\theta(q^{\alpha}z_1\cdots z_n b_1\cdots b_n)
\prod_{1\le i<j\le n}(qz_i/z_j)_\8(qz_j/z_i)_\8}\,\frac{d_qz_1}{z_1}\wedge\cdots\wedge\frac{d_qz_n}{z_n}
=C,
\end{equation*}
which is rewritten as (\ref{eq:01cal I(x)4}) using (\ref{eq:01C}) under the condition (\ref{eq:01beta}). $\square$\\

As a corollary, it is confirmed that the following identity for a contour integral is equivalent to  
the formula (\ref{eq:01cal I(x)4}) of the special case $x=a\in (\mathbb{C}^*)^n$. 
\begin{cor}
Let $\mathbb{T}^n$ be the the direct product of the unit circle, i.e.,  
$
\mathbb{T}^n:=\{(z_1,\ldots,z_n)\in\mathbb{C}^n \,;\,|z_i|=1\}
$. 
Suppose that $|a_i|<1$, $|b_i|<1$ $(i=1,\ldots, n)$ and $a_1\cdots a_n b_1\cdots b_nq^{\alpha+\beta}=q$. 
Then
\begin{eqnarray*}
&&\!\!\!\!\!\!
\Big(\frac{1}{2\pi\sqrt{-1}}\Big)^{\!\!n}\frac{1}{n!}
\int_{\mathbb{T}^n}
\frac{(q^\alpha a_1\cdots a_n z_1^{-1}\cdots z_n^{-1})_\8(q^\beta b_1\cdots b_n z_1\cdots z_n)_\8}
{\prod_{i=1}^n \prod_{j=1}^n (a_iz_j^{-1})_\8(b_i z_j)_\8}\nonumber\\
&&\hskip 80pt 
\times\prod_{1\le i<j\le n}(z_i/z_j)_\8(z_j/z_i)_\8
\frac{dz_1}{z_1}\cdots\frac{dz_n}{z_n}
=\frac{(q^{1-\alpha})_\8(q^{1-\beta})_\8}{(q)_\8^n 
\prod_{i=1}^{n}\prod_{j=1}^{n}(a_ib_j)_\8}.
\end{eqnarray*}
\end{cor}
{\bf Proof.} By residue calculation using (\ref{eq:01cal I(x)4}) in the case $x=a\in (\mathbb{C}^*)^n$. $\square$

%

%
\subsection{Dual expression of the Jackson integral $I(x)$}
For an arbitrary $x=(x_1,x_2,\ldots,x_n)\in (\mathbb{C}^*)^n$ 
we specify $x^{-1}$ as 
\begin{equation}
\label{eq:01x-1}
x^{-1}:=(x_1^{-1},x_2^{-1},\ldots,x_n^{-1})\in (\mathbb{C}^*)^n.
\end{equation}
For the point 
$b=(b_{1},b_{2},\ldots,b_{n})\in (\mathbb{C}^*)^{n}$, 
if we set $y=b^{-1}$ in the connection formula (\ref{eq:01cal I(x)I(y)}), 
then we obtain the expression 
\begin{equation}
\label{eq:01cal I(x)I(b)}
{\cal I}(x)=\frac{\theta(q^\alpha x_1x_2\cdots x_nb_1b_2\cdots b_n)}{\theta(q^\alpha)}{\cal I}(b^{-1}). 
\end{equation}
Since $x=b^{-1}$ is a pole of the function $I(x)$ by definition,  
$I(b^{-1})$ no longer makes sense. 
However, the regularization ${\cal I}(b^{-1})$ appearing on the right-hand side of (\ref{eq:01cal I(x)I(b)}) still has 
meaning as a special value of a holomorphic function. 
We will show a way to realize the regularization ${\cal I}(b^{-1})$ as a computable object by another Jackson integral.  
For this purpose, let $\bar\Phi(z)$ and $\bar \Delta(z)$ be the functions specified by 
\begin{equation}
\label{eq:01barPhi}
\bar\Phi(z):=(z_1z_2\cdots z_n)^{1-\alpha_1-\cdots-\alpha_n-\beta_1-\cdots-\beta_n-\alpha}
\prod_{i=1}^n\prod_{j=1}^n\frac{(qb_j^{-1}z_i)_\8}{(a_jz_i)_\8},
\end{equation}
where $\alpha_i$ and $\beta_i$ are given by $a_i=q^{\alpha_i}, b_i=q^{\beta_i}$, and 
\begin{equation}
\label{eq:01barDelta}
\bar\Delta(z):=\prod_{1\le i<j\le n}(z_i-z_j).
\end{equation}
For $x=(x_1,x_2,\ldots,x_n)\in (\mathbb{C}^*)^n$, we define the sum $\bar I(x)$ by  
\begin{equation}
\label{eq:01bar I(x)}
\bar I(x):=\int_0^{\mbox{\small $x$}\8}
\bar\Phi(z)\bar\Delta(z)\,\frac{d_qz_1}{z_1}\wedge\cdots\wedge\frac{d_qz_n}{z_n},
\end{equation}
which converges under the condition (\ref{eq:01condition01}).
We call $\bar I(x)$ the {\it dual Jackson integral of $I(x)$}, 
and call $\bar I(b)$ its {\it truncation}. 
(The sum $I(x)$ transforms to its dual $\bar I(x)$ up to sign if we interchange the parameters as 
\begin{equation}
\label{eq:01a<-->b}
\alpha\leftrightarrow\beta 
\quad\mbox{and}\quad a_i\leftrightarrow b_i\ (i=1,\ldots,n),
\end{equation}
where $\beta$ is specified by (\ref{eq:01beta}).)
We also define the {\it regularization} $\bar {\cal I}(x)$ of $\bar I(x)$ as 
\begin{equation}
\label{eq:01barh(x)}
\bar {\cal I}(x):=
\frac{\bar I(x)}{\bar h(x)}, \ \mbox{where}\ 
\bar h(x)=(x_1x_2\cdots x_n)^{1-\alpha_1-\cdots-\alpha_{n}-\beta_1-\cdots-\beta_{n}-\alpha}
\frac{\prod_{1\le i<j\le n}x_i\theta(x_j/x_i)}{\prod_{i=1}^n\prod_{j=1}^n\theta(a_jx_i)}.
\end{equation}
In the same manner as Lemma \ref{lem:h-and-s}, we can confirm that the function $\bar {\cal I}(x)$ is also holomorphic and symmetric.

\begin{lem}[reflective equation]
\label{lem:01ref}
The connection between $I(x)$ and $\bar I(x)$ is 
\begin{equation}
\label{eq:01ref1}
I(x)=\frac{h(x)}{\bar h(x^{-1})}\bar I(x^{-1}),\quad 
\mbox{where}\quad
\frac{h(x)}{\bar h(x^{-1})}=
\prod_{i=1}^n\prod_{j=1}^n
x_i^{1-\alpha_j-\beta_j}
\frac{\theta(qa_j^{-1}x_i)}{\theta(b_jx_i)},
\end{equation}
where $x^{-1}$ is specified as in {\rm (\ref{eq:01x-1})}. 
In other words, the relation between ${\cal I}(x)$ and $\bar {\cal I}(x)$ is 
\begin{equation}
\label{eq:01ref2}
{\cal I}(x)=\bar{\cal I}(x^{-1}).
\end{equation}
\end{lem}
{\bf Proof.} 
From the definitions (\ref{eq:01cal I(x)1}) and (\ref{eq:01barh(x)}) the ratio $h(x)/\bar h(x^{-1})$ 
is written as in (\ref{eq:01ref1}). 
Since $
\Delta(z)=(z_1z_2\cdots z_n)^{n-1}\bar\Delta(z^{-1})
$, 
from (\ref{eq:01Phi1}), (\ref{eq:01Delta1}), (\ref{eq:01barPhi}) and (\ref{eq:01barDelta}), 
we have 
\begin{equation}
\label{eq:01PD=hhPD}
\Phi(z)\Delta(z)=
\frac{h(z)}{\bar h(z^{-1})}
\bar\Phi(z^{-1})\bar\Delta(z^{-1}).
\end{equation}
Also since $h(z)/\bar h(z^{-1})$ is invariant under the shift $z_i\to qz_i$,   
by the definitions (\ref{eq:01I(x)1}) and (\ref{eq:01bar I(x)}) of $I(x)$ and $\bar I(x)$, 
the connection (\ref{eq:01ref1}) between $I(x)$ and its dual $\bar I(x)$ is derived from (\ref{eq:01PD=hhPD}). $\square$\\

We use $\bar I(\alpha;x)$ instead of $\bar I(x)$ to see the $\alpha$ dependence. 
From (\ref{eq:01ref1}), the recurrence relation for $\bar I(\alpha;x)$ is completely the same as (\ref{eq:01rec1}) of $I(\alpha;x)$. 
\begin{lem}
\label{lem:rec2}
The function $\bar I(\alpha;x)$ also satisfies the recurrence relation {\rm (\ref{eq:01rec1})} of $I(\alpha;x)$, and is rewritten as  
\begin{equation}
\label{eq:01rec2}
{\bar I}(\alpha;x)=\frac{1-q^{1-\alpha}}
{b_1b_2\cdots b_n(1-q^{1-\alpha}a_1^{-1}a_2^{-1}\cdots a_n^{-1}b_1^{-1}b_2^{-1}\cdots b_n^{-1})}
{\bar I}(\alpha-1;x).
\end{equation}
\end{lem}

We saw above that although $I(b^{-1})$ no longer makes sense,  
its regularization ${\cal I}(b^{-1})$ still has meaning as a special value of a holomorphic function, 
and ${\cal I}(b^{-1})$ is evaluated by the dual integral $\bar{\cal I}(b)$ 
via the reflection equation (\ref{eq:01ref2}). Moreover, by definition the regularization $\bar{\cal I}(b)$ itself is calculated 
using $\bar I(b)$, which is then a truncated Jackson integral.  
Though we have already known the value of ${\cal I}(b^{-1})$ 
through connection formula (\ref{eq:01cal I(x)I(a)}) of the case $x=b^{-1}$, 
the point is that 
we can calculate ${\cal I}(b^{-1})$ directly from $\bar I(b)$, 
whose leading term of its asymptotic behavior as $\alpha\to -\8$
is simply computed as follows. 

\begin{cor}
\label{cor:bar I(a-N;b)}
The asymptotic behavior of $\bar I(\alpha-N;b)$ as $N\to +\8$ is written as
\begin{eqnarray}
\label{eq:01asym2}
\bar I(\alpha-N;b)&\sim&
(1-q)^n(b_1\cdots b_n)^{1-\alpha_1-\cdots-\alpha_{n}-\beta_1-\cdots-\beta_{n}-\alpha+N}
\nonumber\\
&&\times \bar\Delta(b)
\prod_{i=1}^n\prod_{j=1}^n\frac{(qb_i/b_j)_\8}{(a_ib_j)_\8}
\quad(N\to +\8).
\end{eqnarray}
Moreover, by repeated use of {\rm (\ref{eq:01rec2})}, the truncated Jackson integral $\bar I(b)$ is written as  
\begin{equation}
\label{eq:01I(b)}
\bar I(b)=
\frac{
(1-q)^n(b_1\cdots b_n)^
{1-\alpha_1-\cdots-\alpha_n-\beta_1-\cdots-\beta_n-\alpha}
(q)_\8^n 
(q^{1-\alpha})_\8
}
{(q^{1-\alpha}a_1^{-1}\cdots a_n^{-1}b_1^{-1}\cdots b_n^{-1}
)_\8\prod_{i=1}^{n}\prod_{j=1}^{n}(a_ib_j)_\8}
\prod_{1\le i<j\le n}b_i\theta(b_j/b_i).
\end{equation}
\end{cor}
{\bf Proof.}
Using Lemma \ref{lem:rec2}, the arguments are completely parallel to 
Lemma \ref{lem:I(a+N;a)} and Corollary \ref{cor:01I(a)}. 
Actually, if we substitute $a_j$, $b_j$ and $\alpha$ in $\Phi(z)$ of (\ref{eq:01Phi1})
by $a_j\to b_j$, $b_j\to a_j$ and $\alpha\to \beta=1-\alpha_1-\cdots-\alpha_n-\beta_1-\cdots-\beta_n-\alpha$, 
respectively, then $\Phi(z)$ transforms to $\bar\Phi(z)$ in (\ref{eq:01barPhi}), 
so that we obtain the same result as Corollary \ref{cor:01I(a)} with these substitutions. $\square$\\

From (\ref{eq:01cal I(x)I(y)}) and  (\ref{eq:01ref2}), for $x, y\in (\mathbb{C}^*)^{n}$ 
we have the connection formula between ${\cal I}(x)$ and $\bar{\cal I}(y)$ as 
$$
{\cal I}(x)=\frac{\theta(q^\alpha x_1x_2\cdots x_nb_1b_2\cdots b_n)}
{\theta(q^\alpha y_1^{-1}y_2^{-1}\cdots y_n^{-1}b_1b_2\cdots b_n)}
\bar{\cal I}(y).
$$
In particular, if $y=b$, then we have 
$$
{\cal I}(x)=\frac{\theta(q^\alpha x_1x_2\cdots x_nb_1b_2\cdots b_n)}{\theta(q^\alpha)}
\bar{\cal I}(b). 
$$
If we switch the symbols from ${\cal I}(x)$ and $\bar{\cal I}(b)$ to $I(x)$ and $\bar I(b)$, respectively, 
then we obtain 
\begin{equation}
\label{eq:01I(x)I(b)}
I(x)=\frac{h(x)\theta(q^\alpha a_1a_2\cdots a_nb_1b_2\cdots b_n)}{{\bar h}(b)\theta(q^\alpha)}
{\bar I}(b).
\end{equation}
We once again obtain the connection formula between a solution $I(x)$ of 
the $q$-difference equation (\ref{eq:01rec1}) and the special solution $\bar I(b)$ fixed by its 
asymptotic behavior (\ref{eq:01asym2}) as $\alpha\to -\8$, 
as a counterpart of the formula (\ref{eq:01I(x)I(a)}) of the case $\alpha\to +\8$.

The connection formula (\ref{eq:01I(x)I(b)}) with (\ref{eq:01I(b)}) is also another expression 
for the Milne--Gustafson sum, like the formula (\ref{eq:01I(x)I(a)}). 
\\

\noindent
{\bf Remark.} As we have seen above, we used the integrand $\bar\Phi(z)$ instead of $\Phi(z)$, 
which coincides with $\bar\Phi(z)$
up to the $q$-periodic factor $h(z)/h(z^{-1})$, and used the set 
$\{(b_1q^{\nu_1},\ldots,b_nq^{\nu_n})\in (\mathbb{C}^*)^n\,;\,(\nu_1,\ldots,\nu_n)\in \mathbb{N}^n\}$
as the ``$(-\alpha)$-stable cycle" for the dual integral $\bar I(x)$  
when we construct a special solution $\bar I(b)$ expressed by 
(Jackson) integral representation for the $q$-difference equation (\ref{eq:01rec1}) as $\alpha\to -\8$ . 
In the classical setting, this process is usually done by taking an imaginary cycle without changing the integrand $\Phi(z)$
under the ordinary integral representation. In the $q$-analog setting  Aomoto and Aomoto--Kato 
used the integral representation without changing the integrand $\Phi(z)$, but instead, 
they adopted the residue sum on the set 
$\{(b_1^{-1}q^{-\nu_1},\ldots,b_n^{-1}q^{-\nu_n})\in (\mathbb{C}^*)^n\,;\,(\nu_1,\ldots,\nu_n)\in \mathbb{N}^n\}$
of poles of $I(x)$. They call this cycle the {\it $\alpha$-unstable cycle} \cite{Ao94,AK94-2} of $I(x)$ for the parameter $\alpha$.
To carry out this process is called the {\it regularization} in their original paper \cite{Ao90}.
We hope our slight changes of terminology does not bring confusion to the reader. \\

\section{Jackson integral of Dixon--Anderson type}
\label{section:02}
In this section we use $q$-difference equations to show a proof of Evans's summation formula for $q$-Dixon--Anderson integral
introducing a multi-dimensional bilateral extension of his sum,  
which we call the Jackson integral of Dixon--Anderson type. 
In addition, as a limiting case, the Milne--Gustafson summation formula is deduced from it. 
In this sense, the Jackson integral of Dixon--Anderson type as a $q$-analog of the integral (\ref{eq:DA}) 
can also be regarded as a natural multi-variable extension of Ramanujan's $_1\psi_1$ sum .

\subsection{Definitions and the results}
Let $a_1,a_2,\ldots,a_{n+1}$, $b_1,b_2,\ldots,b_{n+1}$ be complex numbers satisfying 
\begin{equation}
\label{eq:02condition01}
q<|a_1a_2\cdots a_{n+1}b_1b_2\cdots b_{n+1}|.
\end{equation}
Let $\Phi(z)$ be specified by 
\begin{equation}
\label{eq:02Phi1}
\Phi(z):=z_1z_2\cdots z_n\prod_{i=1}^n\prod_{j=1}^{n+1}\frac{(qa_j^{-1}z_i)_\8}{(b_jz_i)_\8},
\end{equation}
(cf. (\ref{eq:01Phi1})) and let $\Delta(z)$ be specified by (\ref{eq:01Delta1}).
For $x=(x_1,x_2,\ldots,x_n)\in (\mathbb{C}^*)^n$, we define the sum $J(x)$ by  
\begin{equation}
\label{eq:02J(x)1}
J(x):=\int_0^{\mbox{\small $x$}\8}\Phi(z)\Delta(z)\,\frac{d_qz_1}{z_1}\wedge\cdots\wedge\frac{d_qz_n}{z_n},
\end{equation}
which converges absolutely under the condition (\ref{eq:02condition01}) 
(cf.~\cite[Lemma 3.19]{Gu87}, \cite[Lemma 2.5]{Mi86}). 
We call $J(x)$ the {\it Jackson integral of Dixon--Anderson type}.
By definition $J(x)$ is skew-symmetric. 
\vskip 4mm
\par
For an arbitrary point $(x_1,x_2,\ldots,x_{n+1})\in (\mathbb{C}^*)^{n+1}$, we set the point $(\widehat{x}_i)$ of $(\mathbb{C}^*)^{n}$ by  
$$(\widehat{x}_i):=(x_1,\ldots,x_{i-1},x_{i+1},\ldots,x_{n+1})\in (\mathbb{C}^*)^{n}\quad\mbox{for}\quad i=1,2,\ldots,n+1.$$
%
%
The main result of this section is the following:
\begin{thm}\label{thm:02}
For the function $J(x)$ on $(\mathbb{C}^*)^{n}$ as a sum over $\mathbb{Z}^n$, 
\begin{equation}
\label{eq:02J(x)}
\sum_{i=1}^{n+1}(-1)^{i-1}J(\widehat{x}_i)=
C_0\,\frac{\theta(x_1x_2\cdots x_{n+1}b_1b_2\cdots b_{n+1})}
{\prod_{i=1}^{n+1}\prod_{j=1}^{n+1}\theta(x_ib_j)}
\prod_{1\le i<j\le n+1}x_j\theta(x_i/x_j),
\end{equation}
where $C_0$ is a constant independent of $(x_1,x_2,\ldots,x_{n+1})\in (\mathbb{C}^*)^{n+1}$, which is explicitly written as 
\begin{equation}
\label{eq:02C}
C_0=(1-q)^n
\frac{(q)_\8^n \prod_{i=1}^{n+1}\prod_{j=1}^{n+1}(qa_i^{-1}b_j^{-1})_\8}
{(qa_1^{-1}\cdots a_{n+1}^{-1}b_1^{-1}\cdots b_{n+1}^{-1})_\8}.
\end{equation}
\end{thm}
{\bf Proof.} 
From the denominator of integrand (\ref{eq:02Phi1}), we see the poles of $J(x)$ are included in the set of zero points of 
$\prod_{i=1}^{n+1}\prod_{j=1}^{n+1}\theta(x_ib_j)$. 
Since the sum $\sum_{i=1}^{n+1}(-1)^{i-1}J(\widehat{x}_i)$ is skew-symmetric 
with respect to the permutation of $\{x_1,x_2,\ldots, x_{n+1}\}$, it is divisible 
by $\prod_{1\le i<j\le n+1}x_j\theta(x_i/x_j)$. 
We therefore obtain   
\begin{equation}
\label{eq:02J(x)02}
\sum_{i=1}^{n+1}(-1)^{i-1}J(\widehat{x}_i)=f(x) 
\frac{\prod_{1\le i<j\le n+1}x_j\theta(x_i/x_j)}
{\prod_{i=1}^{n+1}\prod_{j=1}^{n+1}\theta(x_ib_j)},
\end{equation}
where $f(x)=f(x_1,x_2,\ldots,x_{n+1})$ is some holomorphic function on $(\mathbb{C}^*)^{n+1}$. 
Since the left-hand side of (\ref{eq:02J(x)02}) is an invariant under the $q$-shift $x_i\to qx_i$, 
taking account of the quasi-periodicity (\ref{eq:00quasi-period}) of the function on the right-hand of (\ref{eq:02J(x)02}), 
we have that the holomorphic function $f(x)$ must satisfy 
$$T_{x_i}f(z)=-\frac{f(x)}{x_1x_2\cdots x_{n+1}b_1b_2\cdots b_{n+1}}
\quad\mbox{for}\quad i=1,2,\ldots,n+1.$$
This equation has the unique holomorphic solution up to constant, and is written as   
$$f(x)=C_0\,\theta(x_1x_2\cdots x_{n+1} b_1b_2\cdots b_{n+1}),$$
where $C_0$ is some constant independent of $x_1,x_2,\ldots,x_{n+1}$. 
Therefore we obtain the expression (\ref{eq:02J(x)}).
The explicit evaluation of the constant $C_0$ in (\ref{eq:02C}) will be given in Subsection \ref{subsection:02-5}. $\square$ 
\vskip 4mm
\par

We call $J(\widehat{a}_i)$ ($i=1,2,\ldots,n+1$) 
the {\it truncated} Jackson integral, which is defined as a sum over $\mathbb{N}^n$.
As a special case of Theorem \ref{thm:02}, we immediately have the following:
\begin{cor}[Evans \cite{Ev92}]
\label{cor:02}
The {\it truncated} Jackson integral $J(\widehat{a}_i)$ satisfies
\begin{equation}
\label{eq:02J(a)}
\sum_{i=1}^{n+1}(-1)^{i-1}J(\widehat{a}_i)=
(1-q)^n\frac{(q)_\8^n 
(a_1\cdots a_{n+1}b_1\cdots b_{n+1})_\8
}
{\prod_{i=1}^{n+1}\prod_{j=1}^{n+1}(a_ib_j)_\8}
\prod_{1\le i<j\le n+1}a_j\theta(a_i/a_j).
\end{equation}
\end{cor}
{\bf Remark.} 
If we substitute $a_j$ and $b_j$ as 
$a_j\to x_{j-1}$ and $b_j\to q^{s_{j-1}}/x_{j-1}$, respectively, 
then (\ref{eq:02J(a)}) is directly rewritten by 
the ordinary iterated Jackson integral (\ref{eq:00jac1}) according to
\begin{eqnarray}
\label{eq:02evans2}
&&\int_{z_n=x_{n-1}}^{x_{n}}\cdots\int_{z_2=x_1}^{x_2}\int_{z_1=x_0}^{x_1}
\prod_{i=1}^n\prod_{j=0}^{n}\frac{(qz_i/x_j)_\8}{(q^{s_j}z_i/x_j)_\8}\Delta(z)
\,d_qz_1d_qz_2\cdots d_qz_n\nonumber\\
&&=(1-q)^n\frac{(q)_\8^n(q^{s_0+s_1+\cdots+s_{n}})_\8}
{(q^{s_0})_\8(q^{s_1})_\8\cdots (q^{s_{n}})_\8}
\prod_{0\le i<j\le n}\frac{x_j\theta(x_i/x_j)}{(x_iq^{s_j}/x_j)_\8(x_jq^{s_i}/x_i)_\8},
\end{eqnarray}
which exactly coincides with the formula (\ref{eq:Evans1}) established by Evans.
Since (\ref{eq:02evans2}) is already proved in \cite{Ev92}, 
logically speaking, the constant $C_0$ of (\ref{eq:02J(x)}) in Theorem \ref{thm:02} can conversely be evaluated as (\ref{eq:02C}) 
if we use (\ref{eq:02J(a)}), which is equivalent to (\ref{eq:02evans2}),
after putting $(z_1,\ldots,z_{n+1})=(a_1,\ldots,a_{n+1})$ on the equation (\ref{eq:02J(x)}). 
The argument in \cite{Ev92} to prove (\ref{eq:02evans2}) was done under the restriction on $s_i$ as positive integers. 
As is already pointed out in \cite{Ev92}, by analytic continuation $s_i$ can be considered as complex numbers 
after proving (\ref{eq:02evans2}) in the setting $s_i$ are integers. 
But here we prefer to start without restrictions on the parameters being integers, and 
we would like to choose another way for the evaluation of $C_0$
in keeping with our viewpoint. 
Our method then is based on regarding $J(x)$ as a solution of $q$-difference equations fixed by its asymptotic behavior.
Thus, in this paper, (\ref{eq:02J(a)}) is consequently obtained as a corollary via Theorem \ref{thm:02}. \\[4pt]

The remaining subsections are mainly devoted to the evaluation of the constant $C_0$ as (\ref{eq:02C}) in Theorem \ref{thm:02},
which we will see in Subsection \ref{subsection:02-5}.  
In addition, considering the dual expression $\bar J(x)$ of the integral $J(x)$, we can regard 
the Jackson integral of Milne--Gustafson type as a limiting case of that of Dixon--Anderson type. 
(Consequently, we have a proof of the Milne--Gustafson sum again using 
the Jackson integral of Dixon--Anderson type, see Corollary \ref{cor:04}.)
In this sense, the Jackson integral of Dixon--Anderson type as a $q$-analog of the integral (\ref{eq:DA}) 
can also be regarded as a natural multi-dimensional extension of Ramanujan's $_1\psi_1$ sum . \\

For the above purposes we first 
state the $q$-difference equations for 
the regularization and the dual expression of $J(x)$.
\subsection{$q$-difference equations}
Let $\bar\Phi(z)$ be specified by
\begin{equation}
\label{eq:02Phi2}
\bar\Phi(z):=(z_1z_2\cdots z_n)^{1-\alpha_1-\cdots-\alpha_{n+1}-\beta_1-\cdots-\beta_{n+1}}\prod_{i=1}^n\prod_{j=1}^{n+1}\frac{(qb_j^{-1}z_i)_\8}{(a_jz_i)_\8},
\end{equation}
(cf.~(\ref{eq:01barPhi})) where $\alpha_i$ and $\beta_i$ are given by $a_i=q^{\alpha_i}, b_i=q^{\beta_i}$, and let $\bar\Delta(z)$ be specified by
(\ref{eq:01barDelta}).
For $x=(x_1,x_2,\ldots,x_n)\in (\mathbb{C}^*)^n$, we define the sum $\bar J(x)$ by  
\begin{equation}
\label{eq:02bar J(x)}
\bar J(x):=\int_0^{\mbox{\small $x$}\8}
\bar\Phi(z)\bar\Delta(z)\,\frac{d_qz_1}{z_1}\wedge\cdots\wedge\frac{d_qz_n}{z_n},
\end{equation}
which converges under the condition (\ref{eq:02condition01}). 
We call $\bar J(x)$ the {\it dual Jackson integral of} $J(x)$. 
By definition $\bar J(x)$ is skew-symmetric. 
For the specific points 
$$x=(\,\widehat{b}_i)=(b_1,\ldots,b_{i-1},b_{i+1},\ldots,b_{n+1})\in (\mathbb{C}^*)^{n},
\quad i=1,2,\ldots, n+1,$$
we call $\bar J(\,\widehat{b}_i)$ the {\it truncated} Jackson integral, 
which is defined as a sum over $\mathbb{N}^n$.

Let ${\cal J}(x)$ and $h(x)$ be the functions defined by 
\begin{equation}
\label{eq:02h(z)}
{\cal J}(x):=\frac{J(x)}{h(x)},\quad\mbox{where}\quad 
h(x)=x_1x_2\cdots x_n\frac
{\prod_{1\le i<j\le n}x_j\theta(x_i/x_j)}
{\prod_{i=1}^n\prod_{j=1}^{n+1}\theta(b_jx_i)},
\end{equation}
which we call the {\it regularization of $J(x)$}. 
We also define the regularization $\bar {\cal J}(x)$ of $\bar J(x)$ according to 
\begin{equation}
\label{eq:02barh(z)}
\bar {\cal J}(x):=
\frac{\bar J(x)}{\bar h(x)}, \ \mbox{where}\ 
\bar h(x)=(x_1x_2\cdots x_n)^{1-\alpha_1-\cdots-\alpha_{n+1}-\beta_1-\cdots-\beta_{n+1}}
\frac{\prod_{1\le i<j\le n}x_i\theta(x_j/x_i)}{\prod_{i=1}^n\prod_{j=1}^{n+1}\theta(a_jx_i)}.
\end{equation}
Since the trivial poles and zeros of $J(x)$ are canceled out by multiplying together $1/h(x)$ and $J(x)$, 
the function ${\cal J}(x)$ is holomorphic on $(\mathbb{C}^*)^n$, and ${\cal J}(x)$ is symmetric. 
In the same manner, the function $\bar {\cal J}(x)$ is also holomorphic and symmetric.
\begin{lem}[reflective equation]
The connection between $J(x)$ and $\bar J(x)$ is 
\begin{equation}
\label{eq:02ref1}
J(x)=\frac{h(x)}{\bar h(x^{-1})}\bar J(x^{-1}),\quad 
\mbox{where}\quad
\frac{h(x)}{\bar h(x^{-1})}=
\prod_{i=1}^n\prod_{j=1}^{n+1}
x_i^{1-\alpha_j-\beta_j}
\frac{\theta(qa_j^{-1}x_i)}{\theta(b_jx_i)},
\end{equation}
where $x^{-1}$ is specified as in {\rm (\ref{eq:01x-1})}. 
In other words, 
the relation between ${\cal J}(x)$ and $\bar {\cal J}(x)$ is 
\begin{equation}
\label{eq:02ref2}
{\cal J}(x)=\bar{\cal J}(x^{-1}).
\end{equation}
\end{lem}
{\bf Proof.} The proof is the same as that of Lemma \ref{lem:01ref} in Section \ref{section:01}. $\square$\\

We now state the $q$-difference equations for $\bar J(x)$ or $\bar{\cal J}(x)$ under the setting 
$x=(\,\widehat{b}_i)$, $i=1,2,\ldots, n+1$.

\begin{prop}
\label{lem:02q-diff01} 
Fix $x=(\,\widehat{b}_i), i=1,2,\ldots, n+1$ 
and assume 
$1<|a_1a_2\cdots a_{n+1}b_1b_2\cdots b_{n+1}|$. Then 
the recurrence relations for $\bar J(x)$ are given by 
\begin{equation}
\label{eq:02rec2}
T_{a_j}\bar J(x)=(-a_j)^n\frac{\prod_{i=1}^{n+1}(1-b_i^{-1}a_j^{-1})}{1-\prod_{i=1}^{n+1}b_i^{-1}a_i^{-1}}\bar J(x)
\ \ \mbox{and}\ \ 
T_{b_j}\bar J(x)=(-b_j^{-1})^{n}\frac{\prod_{i=1}^{n+1}(1-a_ib_j)}{1-\prod_{i=1}^{n+1}a_ib_i}\bar J(x)
\end{equation}
for $j=1,2,\ldots, n+1.$
The recurrence relations for $\bar {\cal J}(x)$ are given by 
\begin{equation}
T_{a_j}\bar {\cal J}(x)=\frac{\prod_{i=1}^{n+1}(1-b_i^{-1}a_j^{-1})}{1-\prod_{i=1}^{n+1}a_i^{-1}b_i^{-1}}\bar {\cal J}(x)
\quad\mbox{and}\quad
T_{b_j}\bar {\cal J}(x)=\frac{\prod_{i=1}^{n+1}(1-a_i^{-1}b_j^{-1})}{1-\prod_{i=1}^{n+1}a_i^{-1}b_i^{-1}}\bar {\cal J}(x).
\label{eq:02rec-bJ(b)}
\end{equation}
for $j=1,2,\ldots, n+1.$
\end{prop}

The remaining this subsection is devoted to the proof of the above proposition. 
Let $e(c;z)$ be the symmetric polynomial of degree $n$ defined by 
\begin{equation}
\label{eq:02e(c;z)}
e(c;z):=\prod_{i=1}^{n}(1-c^{-1}z_i),
\end{equation}
which has the property that $e(c;z)$ vanishes for $z_i=c$. 
\begin{lem}
\label{lem:q-diff01}
Suppose that $\nabla_{\!i}$ in {\rm(\ref{eq:00nabla})} is defined in terms $\bar\Phi$ instead of $\Phi$. 
If we put $\varphi(z)$ as 
\begin{equation}
\label{eq:02varphi0}
\varphi(z)= z_1^{-1}\prod_{i=1}^{n+1}(b_i-z_1)\times\prod_{1\le j<k\le n}(z_k-b_j),
\end{equation}
then ${\cal A}\nabla_{\!1}\varphi(z)$ is expanded as 
\begin{equation}
\label{eq:02A-nabla-phi0}
{\cal A}\nabla_{\!1}\varphi(z)=\Big(c_0+c_1\frac{e(b_1;z)}{z_1z_2\cdots z_n}\Big)\Delta(z),
\end{equation}
where the constants $c_0$ and $c_1$ are given by 
\begin{equation}
\label{eq:02c0c1}
c_0=-b_1^{-1} a_1^{-1}\cdots a_{n+1}^{-1}\prod_{i=1}^{n+1}(1- a_ib_1),
\quad c_1=(-b_1)^{n}a_1^{-1}\cdots a_{n+1}^{-1}(1-\prod_{i=1}^{n+1}a_ib_i).
\end{equation}
\end{lem}
{\bf Proof.} From (\ref{eq:02Phi2}), the ratio ${T_{z_1}\bar\Phi(z)}/\bar\Phi(z)$ is written as 
$$
\frac{T_{z_1}\bar\Phi(z)}{\bar\Phi(z)}=q\prod_{j=1}^{n+1}\frac{a_j^{-1}-z_1}{b_j-q z_1}.
$$
From (\ref{eq:00nabla}) modified as specified above, 
for the function $\varphi(z)$ defined in (\ref{eq:02varphi0}) we have 
$$
\nabla_{\!1}\varphi(z)=
z_1^{-1}\Big(
\prod_{i=1}^{n+1}(b_i- z_1)
-\prod_{i=1}^{n+1}(a_i^{-1}-z_1)\Big)
\prod_{1\le j<k\le n}(z_k-b_j)
=\frac{\varphi'(z)}{z_1z_2\cdots z_n}
$$
where
\begin{eqnarray}
\label{eq:02nabla-phi0}
\varphi'(z)=\Big(
\prod_{i=1}^{n+1}(b_i- z_1)
-\prod_{i=1}^{n+1}(a_i^{-1}-z_1)\Big)
z_2z_3\cdots z_n\prod_{1\le j<k\le n}(z_k-b_j).
\end{eqnarray}
Since 
$
{\cal A}\nabla_{\!1}\varphi(z)={{\cal A}\varphi'(z)}/{z_1z_2\cdots z_n}
$,
for the purpose of proving (\ref{eq:02A-nabla-phi0}) it suffices to show that 
\begin{equation}
\label{eq:02A-nabla-phi0.5}
{\cal A}\varphi'(z)=\big(c_0z_1z_2\cdots z_n+c_1e(b_1;z)\big)\Delta(z).
\end{equation}
Taking account of the degree of $\varphi'(z)$ as a polynomial of $z$, 
we can expand the skew-symmetrization ${\cal A}\varphi'(z)$ as 
\begin{equation}
\label{eq:02A-nabla-phi1}
{\cal A}\varphi'(z)=\Big(c_0z_1z_2\cdots z_n+\sum_{i=1}^nc_ie(b_i;z)\Big)\Delta(z),
\end{equation}
where the coefficients $c_i$ $(i=0,1,\ldots,n)$ are some constants. 
Here we will confirm that $c_i$ vanishes if $i\ge 2$, and $c_0$ and $c_1$ are evaluated as (\ref{eq:02c0c1}).

First we take $z_1=b_1, z_2=b_2,\ldots, z_n=b_n$. 
Then, from (\ref{eq:02nabla-phi0}) and (\ref{eq:02A-nabla-phi1}) with the vanishing property of $e(b_i;z)$, 
we have two ways of expression of 
$
{\cal A}\varphi'(z)$, 
$$
{\cal A}\varphi'(z)=-\prod_{i=1}^{n+1}(a_i^{-1}- b_1)\times b_2b_3\cdots b_n\Delta(z)=c_0b_1b_2\cdots b_n\Delta(z).
$$
Thus we obtain $c_0=-b_1^{-1}\prod_{i=1}^{n+1}(a_i^{-1}-b_1)
=-b_1^{-1} a_1^{-1}\cdots a_{n+1}^{-1}\prod_{i=1}^{n+1}(1- a_ib_1)$. 

Next we take $z_1=b_1$, $z_2=b_2,\ldots, z_{n-1}=b_{n-1}$ and $z_n\not\in\{b_1,b_2,\ldots,b_{n+1}\}$. Then,
from (\ref{eq:02A-nabla-phi1}) and the vanishing property of $e(a_i;z)$, we have
\begin{equation}
\label{eq:02A-nabla-phi2}
{\cal A}\varphi'(z)=\big(c_0b_1b_2\cdots b_{n-1}z_n+c_n e(b_n;z)\big)\Delta(z).
\end{equation}
On the other hand, from (\ref{eq:02nabla-phi0}), we have 
\begin{equation}
\label{eq:02A-nabla-phi3}
{\cal A}\varphi'(z)
=-\prod_{i=1}^{n+1}(a_i^{-1}- b_1)\times b_2b_3\cdots b_{n-1}z_n\Delta(z)=c_0b_1b_2\cdots b_{n-1}z_n\Delta(z).
\end{equation}
Comparing (\ref{eq:02A-nabla-phi2}) and (\ref{eq:02A-nabla-phi3}), we obtain $c_n=0$. In the same manner, 
by symmetry of ${\cal A}\varphi'(z)$ corresponding to
$\varphi'(z)$ in (\ref{eq:02nabla-phi0}), 
we also obtain $c_i=0$ if $i\ge 2$. Thus we obtain the expansion (\ref{eq:02A-nabla-phi0.5}). 

Lastly, we suppose $z_1=0$ and $z_2=b_2,\ldots, z_{n}=b_{n}$ to determine $c_1$ in (\ref{eq:02A-nabla-phi0.5}). 
Then we have 
\begin{equation}
\label{eq:02A-nabla-phi4}
{\cal A}\varphi'(z)=c_1(-b_1^{-1})^nb_2b_3\cdots b_n\prod_{1\le j<k\le n}(b_k-b_j).
\end{equation}
On the other hand, from (\ref{eq:02nabla-phi0}), we have 
\begin{equation}
\label{eq:02A-nabla-phi5}
{\cal A}\varphi'(z)=\Big(
\prod_{i=1}^{n+1}b_i
-\prod_{i=1}^{n+1}a_i^{-1}\Big)
b_2b_3\cdots b_n\prod_{1\le j<k\le n}(b_k-b_j).
\end{equation}
Comparing (\ref{eq:02A-nabla-phi4}) with (\ref{eq:02A-nabla-phi5}), we obtain $c_1$ as is expressed in (\ref{eq:02c0c1}).  
$\square$\\

\noindent
{\bf Proof of Proposition \ref{lem:02q-diff01}.} 
We will prove (\ref{eq:02rec2}) for $T_{b_j}$ first. 
Without loss of generality, it suffices to show that 
\begin{equation}
\label{eq:02rec0}
T_{b_1}{\bar J}(x)=
(-b_1^{-1})^{n}\frac{\prod_{i=1}^{n+1}(1-a_ib_1)}{1-\prod_{i=1}^{n+1}a_ib_i}{\bar J}(x).
\end{equation}
Since we have 
$
T_{b_1}{\bar\Phi}(z)/{\bar\Phi}(z)=
\prod_{i=1}^{n}(1-b_1^{-1}z_i)/z_i=e(b_1;z)/z_1\cdots z_n
$
by definition, $T_{b_1}{\bar J}(x)$ is expressed by 
$$
T_{b_1}{\bar J}(x)=\int_0^{\mbox{\small $x$}\8}\frac{e(b_1;z)}{z_1z_2\cdots z_n}{\bar\Phi}(z){\bar\Delta}(z)\,\frac{d_qz_1}{z_1}\wedge\cdots\wedge\frac{d_qz_n}{z_n}.
$$
Under the condition $x=(\,\widehat{b}_i)$, $i=1,2,\ldots,n+1$, 
the Jackson integral is truncated, i.e., the support of the Jackson integral in $\mathbb{Z}^n$ is restricted to the fan region $\mathbb{N}^n$. 
Now we assume the condition 
$1<|a_1a_2\cdots a_{n+1}b_1b_2\cdots b_{n+1}|$
on parameters. 
Then
the truncated Jackson integral 
$$
\int_0^{\mbox{\small $x$}\8}\varphi(z){\bar\Phi}(z)\,\frac{d_qz_1}{z_1}\wedge\cdots\wedge\frac{d_qz_n}{z_n}
$$
where $\varphi(z)$ is defined by (\ref{eq:02varphi0}) converges absolutely. 
Therefore, applying (\ref{eq:00A}) in Lemma \ref{lem:00nabla=0} to the fact (\ref{eq:02A-nabla-phi0}) in Lemma \ref{lem:q-diff01}, we obtain the relation
$$c_0{\bar J}(x)+c_1T_{b_1}{\bar J}(x)=0,$$
where $c_0$ and $c_1$ are given in (\ref{eq:02c0c1}). This relation coincides with (\ref{eq:02rec0}). 

Next we will show the $q$-difference equation (\ref{eq:02rec2}) for $T_{a_j}$ of the case $j=1$ 
in the same manner as above. 
Since we have 
$
T_{a_1}{\bar\Phi}(z)/{\bar\Phi}(z)=
\prod_{i=1}^{n}(1-a_1z_i)/z_i=e(a_1^{-1};z)/z_1z_2\cdots z_n$,
$T_{a_1}{\bar J}(x)$ is expressed by 
$$
T_{a_1}{\bar J}(x)=\int_0^{\mbox{\small $x$}\8}\frac{e(a_1^{-1};z)}{z_1z_2\cdots z_n}{\bar \Phi}(z){\bar\Delta}(z)\,\frac{d_qz_1}{z_1}\wedge\cdots\wedge\frac{d_qz_n}{z_n}.
$$
Here if we exchange $b_i$ with $a_i^{-1}$ $(i=1,2,\ldots,n+1)$ 
in the above proof of (\ref{eq:02rec0})
including that of Lemma \ref{lem:q-diff01}, 
the way of argument is completely symmetric for this exchange. 
Therefore (\ref{eq:02rec2}) for $T_{a_1}$ is obtained 
exchanging $b_i$ with $a_i^{-1}$ on the coefficient of  (\ref{eq:02rec0}). 

From the expression (\ref{eq:02barh(z)}) of $\bar h(x)$,  
under the condition $x=(\,\widehat{b}_i)$, $i=1,2,\ldots,n+1$, $\bar h(x)$ satisfies 
\begin{equation*}
T_{a_j}\bar h(x)=(-a_j)^n\bar h(x)\quad\mbox{and}\quad
T_{b_j}\bar h(x)=\frac{b_j}{b_1b_2\cdots b_{n+1}}\bar h(x) 
\quad (j=1,2,\ldots, n+1).
\end{equation*}
Since $\bar {\cal J}(x)={\bar J}(x)/{\bar h}(x)$, from the above equations and (\ref{eq:02rec2}), 
we therefore obtain (\ref{eq:02rec-bJ(b)}). $\square$
\subsection{Evaluation of the truncated Jackson integral }
The main result of this subsection is the evaluation of the regularization of the truncated Jackson integral 
using the $q$-difference equations (\ref{eq:02rec-bJ(b)}) in Proposition \ref{lem:02q-diff01} and its asymptotic behavior for the special direction of parameters. 
\begin{thm}
\label{thm:02main03}
For $x=(\,\widehat{b}_i)$, $i=1,2,\ldots, n+1$, 
the regularized Jackson integral $\bar{\cal J}(x)$ is evaluated as 
\begin{equation}
\label{eq:02cal J(b)}
\bar{\cal J}(\,\widehat{b}_i)=
(1-q)^n
\frac{(q)_\8^n \prod_{j=1}^{n+1}\prod_{k=1}^{n+1}(qa_j^{-1}b_k^{-1})_\8}
{(qa_1^{-1}\cdots a_{n+1}^{-1}b_1^{-1}\cdots b_{n+1}^{-1})_\8}.
\end{equation}
\end{thm}
{\bf Proof.} Without loss of generality, 
it suffices to show (\ref{eq:02cal J(b)}) in the case $x=(\,\widehat{b}_{n+1})$, i.e., the case $x=b=(b_1,\ldots,b_n)\in (\mathbb{C}^*)^{n}$. 
We denote by $C$ the right-hand side of (\ref{eq:02cal J(b)}). 
Then it is immediate to confirm that 
$C$ as a function of $a_j$ and $b_j$ satisfies the same $q$-difference equations as 
(\ref{eq:02rec-bJ(b)}) of $\bar {\cal J}(b)$. Therefore 
the ratio $\bar {\cal J}(b)/C$ is invariant under the $q$-shift with respect to $a_j$ and $b_j$. 

Next, for an integer $N$, let $T^N$ be the $q$-shift operator for a special direction defined as 
$$
T^N: a_i\to q^{-nN}a_i\ (i=1,2,\ldots,n+1);\ b_j \to q^{(n+1)N}b_j\ (j=1,2,\ldots,n);\ b_{n+1}\to q^{-nN}b_{n+1}. 
$$
Then, by definition $T^N \bar J(b)$ is written as
\begin{eqnarray*}
&&T^N \bar J(b)
=(1-q)^n\sum_{(\nu_1,\ldots,\nu_n)\in \mathbb{N}^n} 
(b_1b_2\cdots b_nq^{\nu_1+\cdots+\nu_n+n(n+1)N})
^{1-\alpha_1-\cdots-\alpha_{n+1}-\beta_1-\cdots-\beta_{n+1}+nN} 
\nonumber\\
&&\quad
\times\prod_{i=1}^n\Big(
\frac{(b_{n+1}^{-1}b_iq^{1+\nu_i+(2n+1)N})_\8}{(a_{n+1}b_iq^{\nu_i+N})_\8}
\prod_{j=1}^{n}
\frac{(b_j^{-1}b_iq^{1+\nu_i})_\8}{(a_jb_iq^{\nu_i+N})_\8}\Big)
\prod_{1\le i<j\le n}q^{(n+1)N}(b_iq^{\nu_i}-b_jq^{\nu_j}),
\end{eqnarray*}
so that the leading term of the asymptotic behavior of $T^N \bar J(b)$ as $N\to +\8$ is given by
the term corresponding to $(\nu_1,\ldots,\nu_n)=(0,\ldots,0)$ in the above sum, which is 
\begin{eqnarray}
T^N \bar J(b)
&\sim& (1-q)^n(b_1b_2\cdots b_nq^{n(n+1)N})
^{1-\alpha_1-\cdots-\alpha_{n+1}-\beta_1-\cdots-\beta_{n+1}+nN} 
\nonumber\\
&&\times\prod_{i=1}^n\prod_{j=1}^{n}(qb_j^{-1}b_i)_\8
\prod_{1\le i<j\le n}q^{(n+1)N}(b_i-b_j)
\nonumber\\
&=& 
(b_1b_2\cdots b_nq^{n(n+1)N})
^{1-\alpha_1-\cdots-\alpha_{n+1}-\beta_1-\cdots-\beta_{n+1}+nN}
\nonumber\\
&&\times (1-q)^n(q)_\8^n \prod_{1\le i<j\le n}q^{(n+1)N} b_i\theta(b_j/b_i)
\quad\quad(N\to +\8).
\label{eq:02TNJ}
\end{eqnarray}
On the other hand, from (\ref{eq:02barh(z)}), $\bar h(b)C$ is written as 
\begin{eqnarray*}
\bar h(b)C&=&(b_1b_2\cdots b_n)
^{1-\alpha_1-\cdots-\alpha_{n+1}-\beta_1-\cdots-\beta_{n+1}}(1-q)^n(q)_\8^n \prod_{1\le i<j\le n}b_i\theta(b_j/b_i)\\
&&\times
\frac{\prod_{i=1}^{n+1}(qa_i^{-1}b_{n+1}^{-1})_\8}
{(qa_1^{-1}\cdots a_{n+1}^{-1}b_1^{-1}\cdots b_{n+1}^{-1})_\8
\prod_{i=1}^n\prod_{j=1}^{n+1}(a_jb_i)_\8},
\end{eqnarray*}
so that we have 
\begin{eqnarray}
T^N\Big(\bar h(b)C\Big)&=&(b_1b_2\cdots b_nq^{n(n+1)N})
^{1-\alpha_1-\cdots-\alpha_{n+1}-\beta_1-\cdots-\beta_{n+1}+nN}
\nonumber\\
&&
\times
(1-q)^n(q)_\8^n \prod_{1\le i<j\le n}q^{(n+1)N}b_i\theta(b_j/b_i)
\nonumber\\
&&\times
\frac{\prod_{i=1}^{n+1}(a_i^{-1}b_{n+1}^{-1}q^{1+2nN})_\8}
{(a_1^{-1}\cdots a_{n+1}^{-1}b_1^{-1}\cdots b_{n+1}^{-1}q^{1+nN})_\8
\prod_{i=1}^n\prod_{j=1}^{n+1}(a_jb_iq^N)_\8}.
\label{eq:02TNhC}
\end{eqnarray}
As we saw, the ratio $\bar {\cal J}(b)/C$ is invariant under the $q$-shift with respect to $a_j$ and $b_j$.
Thus $\bar {\cal J}(b)/C$ is also invariant under the $q$-shift $T^N$.
Therefore, comparing (\ref{eq:02TNJ}) with (\ref{eq:02TNhC}), we obtain
\begin{eqnarray*}
\frac{\bar {\cal J}(b)}{C}&=&T^N\frac{\bar {\cal J}(b)}{C}=\frac{T^N\bar J(b)}{T^N\bar h (b)C}
=
\lim_{N\to +\8}\frac{T^N\bar J(b)}{T^N\bar h (b)C}\\[2pt]
&=&\lim_{N\to +\8}
\frac
{(a_1^{-1}\cdots a_{n+1}^{-1}b_1^{-1}\cdots b_{n+1}^{-1}q^{1+nN})_\8
\prod_{i=1}^n\prod_{j=1}^{n+1}(a_jb_iq^N)_\8}
{\prod_{i=1}^{n+1}(a_i^{-1}b_{n+1}^{-1}q^{1+2nN})_\8}
\\[2pt]
&=&1,
\end{eqnarray*}
and thus $\bar {\cal J}(b)=C$. $\square$
\subsection{A remark on the relation between $\bar J(b)$ and $\bar I(b)$}
As an application of the $q$-difference equations (\ref{eq:02rec2}) for $\bar J(x)$,  
we can show that 
the product formula (\ref{eq:01I(b)}) for the Milne--Gustafson sum $\bar I(b)$ in Corollary \ref{cor:bar I(a-N;b)}
(or (\ref{eq:01I(a)}) of $I(a)$ in Corollary \ref{cor:01I(a)} by their duality (\ref{eq:01a<-->b}) of parameters) is 
a special case of (\ref{eq:02cal J(b)}) in Theorem \ref{thm:02main03}. 
This indicates a way to prove the Milne--Gustafson summation formula from 
the product formula of the Jackson integral of Dixon--Anderson type. 
\begin{cor}
\label{cor:04} For $b=(b_1,b_2,\ldots,b_n)$, the truncated Jackson integral $\bar J(b)$ of Dixon--Anderson type 
is expressed as 
\begin{equation}
\label{eq:02JandI}
\bar J(b)=\bar I(b)\prod_{i=1}^n\frac{(qa_i^{-1}b_{n+1}^{-1})_\8}{(b_ia_{n+1})_\8},
\end{equation}
where $\bar I(b)$ is the truncated Jackson integral of Milne--Gustafson type 
defined by {\rm (\ref{eq:01bar I(x)})} with the setting 
$\alpha=\alpha_{n+1}+\beta_{n+1}$. 
In particular, $\bar I(b)$ is expressed as {\rm (\ref{eq:01I(b)})} in Corollary \ref{cor:bar I(a-N;b)}.
\end{cor}
{\bf Remark.} 
From (\ref{eq:02JandI}), $\bar I(b)$ is a limiting case of $\bar J(b)$ with the $q$-shift $a_{n+1}\to q^N a_{n+1}$
and $b_{n+1}\to q^{-N} b_{n+1}$ $(N\to +\8)$. 
Conversely, the product formula (\ref{eq:02cal J(b)}) of $\bar J(b)$ in Theorem \ref{thm:02main03} is reconstructed from 
the product formula (\ref{eq:01I(b)}) of $\bar I(b)$ 
via the connection (\ref{eq:02JandI}). \\

\noindent
{\bf Proof.} 
From (\ref{eq:02rec2}) the recurrence relation of $\bar J(b)$ with respect to the $q$-shift $a_{n+1}\to qa_{n+1}$
and $b_{n+1}\to q^{-1} b_{n+1}$ is written as 
\begin{equation*}
\bar J(b)=T_{b_{n+1}}^{-1}T_{a_{n+1}}\bar J(b)\times 
\prod_{i=1}^{n}\frac{1-qa_i^{-1}b_{n+1}^{-1}}{1-b_ia_{n+1}}.
\end{equation*}
By repeated use of this equation we have 
\begin{equation}
\label{eq:02JandI1}
\bar J(b)=T_{b_{n+1}}^{-N}T_{a_{n+1}}^N\bar J(b)\times
\prod_{i=1}^n\frac{(qa_i^{-1}b_{n+1}^{-1})_N}{(b_ia_{n+1})_N}
=\lim_{N\to \8}T_{b_{n+1}}^{-N}T_{a_{n+1}}^N\bar J(b)\times
\prod_{i=1}^n\frac{(qa_i^{-1}b_{n+1}^{-1})_\8}{(b_ia_{n+1})_\8}.
\end{equation}
Moreover, by definition $\displaystyle \lim_{N\to \8}T_{a_{n+1}}^NT_{b_{n+1}}^{-N}\bar J(b)$ is written as 
\begin{eqnarray}
\lim_{N\to \8}T_{b_{n+1}}^{-N}T_{a_{n+1}}^N\bar J(b)
&=&\lim_{N\to \8}
(1-q)^n\!\!\!\!\!
\sum_{(\nu_1,\ldots,\nu_n)\in \mathbb{N}^n}
(b_1b_2\cdots b_n q^{\nu_1+\cdots+\nu_n})^{1-\alpha_1-\cdots-\alpha_{n+1}-\beta_1-\cdots-\beta_{n+1}}
\nonumber\\
&&
\times\prod_{i=1}^n\Big(
\frac{(b_{n+1}^{-1}b_iq^{1+\nu_i+N})_\8}{(a_{n+1}b_iq^{\nu_i+N})_\8}
\prod_{j=1}^{n}\frac{(b_j^{-1}b_iq^{1+\nu_i})_\8}{(a_jb_iq^{\nu_i})_\8}\Big)
\prod_{1\le i<j\le n}(b_iq^{\nu_i}-b_jq^{\nu_j})
\nonumber\\
&=&
(1-q)^n\!\!\!\!\!
\sum_{(\nu_1,\ldots,\nu_n)\in \mathbb{N}^n}
(b_1b_2\cdots b_n q^{\nu_1+\cdots+\nu_n})^{1-\alpha_1-\cdots-\alpha_{n+1}-\beta_1-\cdots-\beta_{n+1}}
\nonumber\\
&&
\times\prod_{i=1}^n
\prod_{j=1}^{n}\frac{(b_j^{-1}b_iq^{1+\nu_i})_\8}{(a_jb_iq^{\nu_i})_\8}
\prod_{1\le i<j\le n}(b_iq^{\nu_i}-b_jq^{\nu_j}),
\label{eq:02JandI2}
\end{eqnarray}
which exactly coincides with the definition of $\bar I(b)$ under the setting $\alpha=\alpha_{n+1}+\beta_{n+1}$. 
From (\ref{eq:02JandI1}) and (\ref{eq:02JandI2}), we therefore obtain (\ref{eq:02JandI}). 

Next, using (\ref{eq:02JandI}), (\ref{eq:02h(z)}) and (\ref{eq:02cal J(b)}) of Theorem \ref{thm:02main03}, 
the sum $\bar I(b)$ is conversely calculated as 
\begin{eqnarray*}
\bar I(b)&=&\prod_{i=1}^n\frac{(b_ia_{n+1})_\8}{(qa_i^{-1}b_{n+1}^{-1})_\8}\bar J(b)
=\prod_{i=1}^n\frac{(b_ia_{n+1})_\8}{(qa_i^{-1}b_{n+1}^{-1})_\8}\bar{\cal J}(b)\bar h(b)\\
&=&
\prod_{i=1}^n\frac{(b_ia_{n+1})_\8}{(qa_i^{-1}b_{n+1}^{-1})_\8}
\times(1-q)^n
\frac{(q)_\8^n \prod_{i=1}^{n+1}\prod_{j=1}^{n+1}(qa_i^{-1}b_j^{-1})_\8}
{(qa_1^{-1}\cdots a_{n+1}^{-1}b_1^{-1}\cdots b_{n+1}^{-1})_\8}\\
&&\quad\times(b_1b_2\cdots b_n)^{1-\alpha_1-\cdots-\alpha_{n+1}-\beta_1-\cdots-\beta_{n+1}}
\frac{\prod_{1\le i<j\le n}b_i\theta(b_j/b_i)}
{\prod_{i=1}^{n}\prod_{j=1}^{n+1}\theta(a_jb_i)}
\\
&=&(1-q)^n 
\frac{(b_1b_2\cdots b_n)^{1-\alpha_1-\cdots-\alpha_{n+1}-\beta_1-\cdots-\beta_{n+1}}}
{(qa_1^{-1}\cdots a_{n+1}^{-1}b_1^{-1}\cdots b_{n+1}^{-1})_\8}
\frac{(q)_\8^n (qa_{n+1}^{-1}b_{n+1}^{-1})_\8}
{\prod_{i=1}^{n}\prod_{j=1}^{n}(a_ib_j)_\8}
\prod_{1\le i<j\le n}b_i\theta(b_j/b_i),
\end{eqnarray*}
which exactly coincides with (\ref{eq:01I(b)}) of 
Corollary \ref{cor:bar I(a-N;b)} under the setting $\alpha=\alpha_{n+1}+\beta_{n+1}$. 
$\square$

\subsection{Evaluation of the constant $C_0$}
\label{subsection:02-5} 
In this subsection we will evaluate the constant $C_0$ in (\ref{eq:02J(x)}) of Theorem \ref{thm:02}.  
From the expression (\ref{eq:02J(x)}), the constant $C_0$ is written as the sum of functions, 
\begin{equation}
\label{eq:02C-sum1}
C_0=\sum_{k=1}^{n+1}(-1)^{k-1}g_k(x_1,x_2,\ldots,x_{n+1}),
\end{equation}
where 
$$\displaystyle 
g_k(x_1,x_2,\ldots,x_{n+1})
:=J(\widehat{x}_k)
\frac
{\prod_{i=1}^{n+1}\prod_{j=1}^{n+1}\theta(x_ib_j)}
{\theta(x_1x_2\cdots x_{n+1}b_1b_2\cdots b_{n+1})\prod_{1\le i<j\le n+1}x_j\theta(x_i/x_j)}
.
$$
Using (\ref{eq:02h(z)}) and (\ref{eq:02ref2}), we have 
\begin{eqnarray*}
g_k(x_1,x_2,\ldots,x_{n+1})&=&
\frac
{\bar {\cal J}(\widehat{x}_k^{\,-1})h(\widehat{x}_k)\prod_{i=1}^{n+1}\prod_{j=1}^{n+1}\theta(x_ib_j)}
{\theta(x_1x_2\cdots x_{n+1}b_1b_2\cdots b_{n+1})\prod_{1\le i<j\le n+1}x_j\theta(x_i/x_j)}\\
&=&
\frac
{\bar {\cal J}(\widehat{x}_k^{\,-1})x_1\cdots x_{k-1}x_{k+1}\cdots x_{n+1}\prod_{j=1}^{n+1}\theta(x_kb_j)}
{\theta(x_1x_2\cdots x_{n+1}b_1b_2\cdots b_{n+1})
\prod_{i=1}^{k-1}x_k\theta(x_i/x_k)\prod_{j=k+1}^{n+1}x_j\theta(x_k/x_j)}\\
&=&
(-1)^{k-1}\bar {\cal J}(\widehat{x}_k^{\,-1})
\frac
{\theta(x_kb_k)}
{\theta(x_1x_2\cdots x_{n+1}b_1b_2\cdots b_{n+1})
}
\prod_{1\le i\le n+1\atop i\ne k}\frac{\theta(x_kb_i)}{\theta(x_k/x_i)}.
\end{eqnarray*}
Using this and (\ref{eq:02C-sum1}) shows
\begin{equation}
\label{eq:02C-sum2}
C_0=\sum_{k=1}^{n+1}
\bar {\cal J}(\widehat{x}_k^{\,-1})
\frac
{\theta(x_kb_k)}
{\theta(x_1x_2\cdots x_{n+1}b_1b_2\cdots b_{n+1})
}
\prod_{1\le i\le n+1\atop i\ne k}\frac{\theta(x_kb_i)}{\theta(x_k/x_i)}.
\end{equation}
Since $C_0$ is a constant independent of $(x_1,x_2,\ldots,x_{n+1})$, 
it suffices to calculate the right-hand side of (\ref{eq:02C-sum2}) in the specific case 
$$
x_1=b_1^{-1},x_2=b_2^{-1},\cdots,x_n=b_n^{-1}
\quad\mbox{and}\quad x_{n+1}\not\in \{b_1^{-1},b_2^{-1},\ldots,b_{n+1}^{-1}\}.
$$
We now impose this condition. Noting $\theta(x_kb_k)=0$ if $k\ne n+1$ we see from (\ref{eq:02C-sum2}) that
$$
C_0=\bar{\cal J}(\,\widehat{b}_{n+1})
\frac
{\theta(x_{n+1}b_{n+1})}
{\theta(b_1^{-1}b_2^{-1}\cdots b_n^{-1}x_{n+1}b_1b_2\cdots b_{n+1})
}
\prod_{i=1}^n\frac{\theta(x_{n+1}b_i)}{\theta(x_{n+1}/b_i^{-1})}
=\bar{\cal J}(\,\widehat{b}_{n+1}).
$$
Since $\bar {\cal J}(\,\widehat{b}_{n+1})$ has already evaluated as (\ref{eq:02cal J(b)}) in Theorem \ref{thm:02main03}, we obtain $C_0$ 
as is expressed in (\ref{eq:02C}) of Theorem \ref{thm:02}. $\square$

\section{Jackson integral of Gustafson's $A_n$-type}
\label{section:03}
In this section we show a proof for the product formula of Gustafson's $A_n$ sum. 
As is pointed out in \cite{Gu87}, the Milne--Gustafson sum $I(x)$ can also be deduced from Gustafson's $A_n$ sum. 
\subsection{Definitions and results}
Let $a_i, b_i$ $(1\le i\le n+1)$ and $d$ be complex numbers in $\mathbb{C}^*$. 
In this section we define $\Phi(z)$ and $\Delta(z)$ by 
\begin{equation}
\label{eq:03Phi0}
\Phi(z):=\prod_{i=1}^{n+1}
\prod_{j=1}^{n+1}\frac{(qa_j^{-1}z_i)_\8}{(b_jz_i)_\8},\quad
\Delta(z):=\prod_{1\le i<j\le n+1}(z_j-z_i),
\end{equation}
under the condition 
\begin{equation}
\label{eq:03balance}
z_1z_2\cdots z_{n+1}=d.
\end{equation}
(The functions $\Phi(z)$ and $\Delta(z)$ are regarded as a function on $(\mathbb{C}^*)^n$ of $n$ variables 
putting $z_{n+1}=dz_1^{-1}\cdots z_n^{-1}$ into the definition (\ref{eq:03Phi0}).) 
For $x=(x_1,x_2,\ldots,x_n)\in (\mathbb{C}^*)^n$, we define the sum $K(x)$ by  
\begin{equation}
\label{eq:03K(x)1}
K(x):=\int_0^{\mbox{\small $x$}\8}\Phi(z)\Delta(z)\,\frac{d_qz_1}{z_1}\wedge\cdots\wedge\frac{d_qz_n}{z_n},
\end{equation}
which converges absolutely under the condition 
\begin{equation*}
\label{eq:03condition01}
q<|a_1a_2\cdots a_{n+1}b_1b_2\cdots b_{n+1}|.
\end{equation*}
(See \cite[Lemma 3.19]{Gu87}.)
We call $K(x)$ the {\it Jackson integral of Gustafson's $A_n$-type}.  
If we set $x=a=(a_1,a_2,\ldots,a_n)$ on (\ref{eq:03K(x)1}), then we call $K(a)$ the {\it truncated} Jackson integral. 
For an arbitrary $x=(x_1,x_2,\ldots,x_n)\in (\mathbb{C}^*)^n$ we define the function $h(x)$ by 
\begin{equation}
\label{eq:03h(x)1}
h(x):=(-1)^n\frac{\prod_{1\le i<j\le n}x_j\theta(x_i/x_j)}
{\prod_{i=1}^{n}\prod_{j=1}^{n+1}\theta(b_jx_i)}
\frac{\prod_{i=1}^{n}x_i\theta(dx_1^{-1}\cdots x_n^{-1}/x_i)}
{\prod_{j=1}^{n+1}\theta(b_jdx_1^{-1}\cdots x_n^{-1})}.
\end{equation}
We remark that, if we denote by $x_{n+1}$ the combination of variables $dx_1^{-1}x_2^{-1}\cdots x_{n}^{-1}$, then $h(x)$ is expressed as  
\begin{equation}
\label{eq:03h(x)2}
h(x)=\frac{\prod_{1\le i<j\le n+1}x_j\theta(x_i/x_j)}
{\prod_{i=1}^{n+1}\prod_{j=1}^{n+1}\theta(b_jx_i)}.
\end{equation}
We define the {\it regularization} ${\cal K}(x)$ of the sum $K(x)$ by ${\cal K}(x):=K(x)/h(x)$. 
Since $h(x)$ is invariant under the $q$-shift $x_i\to qx_i$ $(1\le i\le n)$, 
the regularization ${\cal K}(x)$ is also expressed as 
\begin{equation*}
{\cal K}(x)=\frac{K(x)}{h(x)}
=\int_0^{\mbox{\small $x$}\8}\frac{\Phi(z)\Delta(z)}{h(z)}
\,\frac{d_qz_1}{z_1}\wedge\cdots\wedge\frac{d_qz_n}{z_n},
\end{equation*}
so that, for $x=(x_1,x_2,\ldots,x_n)\in (\mathbb{C}^*)^n$ we have a 
Macdonald-type sum expression (recall (\ref{eq:01cal I(x)4})) of the regularization ${\cal K}(x)$ as 
\begin{equation}
\label{eq:03cal K(x)1}
{\cal K}(x)=\int_0^{\mbox{\small $x$}\8}\frac{\prod_{i=1}^{n+1}\prod_{j=1}^{n+1} (qa_j^{-1}z_i)_\8(qb_j^{-1}z_i^{-1})_\8}
{\prod_{1\le i<j\le n+1}(qz_i/z_j)_\8(qz_j/z_i)_\8}\,\frac{d_qz_1}{z_1}\wedge\cdots\wedge\frac{d_qz_n}{z_n},
\end{equation}
where the variable $z_{n+1}$ in the integrand satisfies the condition (\ref{eq:03balance}) by definition. 
\begin{lem} 
\label{lem:03hol}
The $q$-periodic function 
${\cal K}(x)$ is holomorphic on $(\mathbb{C}^*)^n$, and is consequently a constant independent of $x\in (\mathbb{C}^*)^n$. 
\end{lem}
{\bf Remark.}
As we will see in Lemma \ref{lem:03q-diff1} below, 
${\cal K}(x)$ satisfies the $q$-difference system of rank 1, which is independent of $x\in (\mathbb{C}^*)^n$. 
We regard the connection coefficient between a general solution ${\cal K}(x)$ of the system and a special solution ${\cal K}(a)$ as 1, i.e., 
\begin{equation}
\label{eq:03cK(x)=cK(a)}
{\cal K}(x)={\cal K}(a),
\end{equation}
for an arbitrary $x\in (\mathbb{C}^*)^n$, which is just the statement of Lemma \ref{lem:03hol}. 
In the same way, for $b^{-1}$ as specified in (\ref{eq:01x-1})  we can also consider 
the equation ${\cal K}(x)={\cal K}(b^{-1})$ as a connection formula.\\[2pt]
\noindent
{\bf Proof.} We temporarily write $x_{n+1}=dx_1^{-1}x_2^{-1}\cdots x_n^{-1}$. 
From the expression (\ref{eq:03cal K(x)1}), the function ${\cal K}(x)$ has the poles lying only in the set 
$\{x=(x_1,x_2,\ldots,x_n)\in (\mathbb{C}^*)^n\,;\,\prod_{1\le i<j\le n+1}\theta(x_i/x_j)=0\}$. 
Then ${\cal K}(x)$ is written as 
\begin{equation}
\label{eq:03cal K(x)2}
{\cal K}(x)=\frac{f(x)}{\prod_{1\le i<j\le n+1}x_j\theta(x_i/x_j)}
\end{equation}
where $f(x)$ is some holomorphic function on $(\mathbb{C}^*)^n$.
Since ${\cal K}(x)$ and $\prod_{1\le i<j\le n+1}x_j\theta(x_i/x_j)$ 
are symmetric and skew-symmetric for permutation of $x_1,x_2,\ldots,x_{n+1}$, respectively, 
the holomorphic function $f(x)$ is skew-symmetric. Thus $f(x)$ vanishes if $x_i=x_j$. 
This means that $f(x)$ is divisible by $x_j\theta(x_i/x_j)$. 
Since the holomorphic function $f(x)$ is divisible by $\prod_{1\le i<j\le n+1}x_j\theta(x_i/x_j)$, 
from the expression (\ref{eq:03cal K(x)2}), ${\cal K}(x)$ is holomorphic function on $(\mathbb{C}^*)^n$. 
By definition of Jackson integrals, the holomorphic function ${\cal K}(x)$ is invariant under the $q$-shift $x_i\to qx_i$
$(1\le i\le n)$. Therefore ${\cal K}(x)$ is a constant independent of $x$. $\square$\\

The aim of this section is the evaluation of ${\cal K}(x)$ as a constant.  
The main theorem is stated as follows:
\begin{thm}[Gustafson \cite{Gu87}] 
\label{thm:03main}
For an arbitrary $x=(x_1,x_2,\cdots,x_n)\in (\mathbb{C}^*)^n$, the regularization ${\cal K}(x)$ is a constant independent of $x$, which is expressed as 
\begin{equation}
\label{eq:03cal K(x)3}
{\cal K}(x)
=(1-q)^n
\frac{(q)_\8^n(qa_1^{-1}\cdots a_{n+1}^{-1}d)_\8
(qb_1^{-1}\cdots b_{n+1}^{-1}d^{-1})_\8}
{(qa_1^{-1}\cdots a_{n+1}^{-1}b_1^{-1}\cdots b_{n+1}^{-1})_\8}
\prod_{i=1}^{n+1}\prod_{j=1}^{n+1}(qa_i^{-1}b_j^{-1})_\8.
\end{equation}
In other words, from the expression $K(x)={\cal K}(x)h(x)$, the sum $K(x)$ is expressed as 
\begin{eqnarray}
\label{eq:03K(x)2}
K(x)
&=&(1-q)^n
\frac{(q)_\8^n(qa_1^{-1}\cdots a_{n+1}^{-1}d)_\8
(qb_1^{-1}\cdots b_{n+1}^{-1}d^{-1})_\8}
{(qa_1^{-1}\cdots a_{n+1}^{-1}b_1^{-1}\cdots b_{n+1}^{-1})_\8}
\prod_{i=1}^{n+1}\prod_{j=1}^{n+1}(qa_i^{-1}b_j^{-1})_\8
\nonumber\\
&&\times
(-1)^n\frac{\prod_{1\le i<j\le n}x_j\theta(x_i/x_j)}
{\prod_{i=1}^{n}\prod_{j=1}^{n+1}\theta(b_jx_i)}
\frac{\prod_{i=1}^{n}x_i\theta(dx_1^{-1}\cdots x_n^{-1}/x_i)}
{\prod_{j=1}^{n+1}\theta(b_jdx_1^{-1}\cdots x_n^{-1})}.
\end{eqnarray}
\end{thm}
{\bf Proof.} Using (\ref{eq:03cK(x)=cK(a)}), the evaluation of ${\cal K}(x)$ is deduced from the special case ${\cal K}(a)$ of 
$a=(a_1,a_2,\cdots,a_n)\in (\mathbb{C}^*)^n$, 
which will be shown in Lemma \ref{lem:03main} below. 
We will also mention another way to evaluate ${\cal K}(x)$ in Subsection \ref{subsection:03-3}. 
$\square$\\[5pt]
{\bf Remark 1.} 
If we write $x_{n+1}=dx_1^{-1}x_2^{-1}\cdots x_n^{-1}$, using (\ref{eq:03h(x)2}), the evaluation (\ref{eq:03K(x)2}) of $K(x)$ 
reads more symmetrical
\begin{eqnarray}
\label{eq:03K(c)}
K(x)&=&(1-q)^n
\frac{(q)_\8^n(qa_1^{-1}\cdots a_{n+1}^{-1}x_1\cdots x_{n+1})_\8
(qb_1^{-1}\cdots b_{n+1}^{-1}x_1^{-1}\cdots x_{n+1}^{-1})_\8}
{(qa_1^{-1}\cdots a_{n+1}^{-1}b_1^{-1}\cdots b_{n+1}^{-1})_\8}
\nonumber\\
&&\times
\prod_{i=1}^{n+1}\prod_{j=1}^{n+1}(qa_i^{-1}b_j^{-1})_\8
\frac{\prod_{1\le i<j\le n+1}x_j\theta(x_i/x_j)}
{\prod_{i=1}^{n+1}\prod_{j=1}^{n+1}\theta(b_ix_j)}.
\end{eqnarray}
It is easy to confirm that (\ref{eq:03K(c)}) exactly coincides with Gustafson's original evaluation \cite{Gu87}.
For our context it is very important to clearly distinguish the part independent of $x$ from that dependent on $x$.
This is why we use the expressions (\ref{eq:03cal K(x)3}) or (\ref{eq:03K(x)2}) rather than (\ref{eq:03K(c)}). 
\\[4pt]
{\bf Remark 2.} As is pointed out in Gustafson's original paper \cite[p.1593]{Gu87}, 
the Milne--Gustafson sum $I(x)$ discussed in Section \ref{section:01} 
is immediately deduced from Gustafson's $A_n$ sum $K(x)$.
In this sense, Gustafson's $A_n$ sum is regarded as an extension of the Milne--Gustafson sum.\\[4pt]

Though the following lemma is a special case of Theorem \ref{thm:03main}, from the fact (\ref{eq:03cK(x)=cK(a)}), 
logically it suffices to evaluate ${\cal K}(a)$ instead of proving Theorem \ref{thm:03main}.
\begin{lem}
\label{lem:03main}
For the special point $x=a=(a_1,a_2,\cdots,a_n)$, the truncation ${\cal K}(a)$ is expressed as 
\begin{equation}
\label{eq:03cal K(a)}
{\cal K}(a)
=(1-q)^n
\frac{(q)_\8^n(qa_1^{-1}\cdots a_{n+1}^{-1}d)_\8
(qb_1^{-1}\cdots b_{n+1}^{-1}d^{-1})_\8}
{(qa_1^{-1}\cdots a_{n+1}^{-1}b_1^{-1}\cdots b_{n+1}^{-1})_\8}
\prod_{i=1}^{n+1}\prod_{j=1}^{n+1}(qa_i^{-1}b_j^{-1})_\8. 
\end{equation}
\end{lem}
{\bf Proof.} Subsection \ref{subsection:03-4} will be devoted to the proof of this lemma, computing 
the asymptotic behavior of the truncated Jackson integral $K(a)$. $\square$
\\[4pt]
\noindent
{\bf Remark.} Notice that 
$a_{n+1}$ fixed in the definition of $\Phi(z)$ does not necessarily satisfy the relation 
$d=a_1a_2\cdots a_{n+1}$. Thus in general, $h(a)$ does not coincide with 
\begin{equation}
\label{eq:03h(a)1}
\frac{\prod_{1\le i<j\le n+1}a_j\theta(a_i/a_j)}
{\prod_{i=1}^{n+1}\prod_{j=1}^{n+1}\theta(a_ib_j)}, 
\end{equation}
and the correct expression of $h(a)$ follows by definition as
\begin{equation}
\label{eq:03h(a)2}
h(a)=(-1)^n
\frac{\prod_{1\le i<j\le n}a_j\theta(a_i/a_j)}
{\prod_{i=1}^{n}\prod_{j=1}^{n+1}\theta(b_ja_i)}
\frac{\prod_{i=1}^{n}a_i\theta(da_1^{-1}\cdots a_n^{-1}/a_i)}
{\prod_{j=1}^{n+1}\theta(b_jda_1^{-1}\cdots a_n^{-1})},
\end{equation}
which will be used later. 
However, if we set $d=a_1a_2\cdots a_{n+1}$, then 
$h(a)$ is written as (\ref{eq:03h(a)1}), 
so that $K(a)={\cal K}(a)h(a)$ is written more simply as 
$$
K(a)
=(1-q)^n
\frac{(q)_\8^{n+1}\prod_{1\le i<j\le n+1}a_j\theta(a_i/a_j)}
{\prod_{i=1}^{n+1}\prod_{j=1}^{n+1}(a_ib_j)_\8},
$$
which is referred as {\it Milne's fundamental theorem of $U(n)$ series} \cite{Mi85} in \cite[(2.2) p.421]{Ro04}.

\subsection{$q$-difference equations}
In this subsection we show the $q$-difference equations for $K(x)$ with respect to parameters explicitly. 
\begin{lem}
\label{lem:03q-diff1}
The recurrence relations for $K(x)$ are given by 
\begin{eqnarray}
T_{a_j}K(x)&=&K(x)
\frac{1-d\prod_{i=1}^{n+1}a_i^{-1}}
{1-\prod_{i=1}^{n+1}a_i^{-1}b_i^{-1}}
\prod_{i=1}^{n+1}(1-b_i^{-1}a_j^{-1}),\label{eq:03rec01}\\
T_{b_j}K(x)&=&K(x)\frac{1-d\prod_{i=1}^{n+1}b_i}{1-\prod_{i=1}^{n+1}a_ib_i}\prod_{i=1}^{n+1}(1-a_ib_j),
\label{eq:03rec02}
\end{eqnarray}
for $j=1,2,\ldots, n+1$. 
The recurrence relations for ${\cal K}(x)$ are given by 
\begin{eqnarray}
T_{a_j}{\cal K}(x)&=&{\cal K}(x)
\frac{1-d\prod_{i=1}^{n+1}a_i^{-1}}
{1-\prod_{i=1}^{n+1}a_i^{-1}b_i^{-1}}
\prod_{i=1}^{n+1}(1-b_i^{-1}a_j^{-1}),
\label{eq:03rec03}\\
T_{b_j}{\cal K}(x)&=&{\cal K}(x)
\frac{1-d^{-1}\prod_{i=1}^{n+1}b_i^{-1}}
         {1-\prod_{i=1}^{n+1}a_i^{-1}b_i^{-1}}
\prod_{i=1}^{n+1}(1-a_i^{-1}b_j^{-1}),
\label{eq:03rec04}
\end{eqnarray}
for $j=1,2,\ldots, n+1.$
\end{lem}

The rest of this subsection is devoted to proving the above $q$-difference equations.
Before proving Lemma \ref{lem:03q-diff1}, we will show two lemmas. 
In this section, $\Phi(z)$ in (\ref{eq:03Phi0}) is considered 
under the restriction (\ref{eq:03balance}). Involving this restriction, 
we need a slight change of the basic lemma (Lemma \ref{lem:00nabla=0}) 
to Lemma \ref{lem:03nabla=0} stated below in order to suit the setting of this section. 

For a function $\varphi(z)$ of $z=(z_1,z_2,\ldots,z_n) \in(\mathbb{C}^*)^n$, 
let $\nabla_{\!i,j}\varphi(z)$ ($1\le i,j\le n+1$) 
be the functions defined by
\begin{equation}
\label{eq:03nabla}
(\nabla_{\!i,j}\varphi)(z):=\varphi(z)-\frac{T_{z_j}^{-1}T_{z_i}\Phi(z)}{\Phi(z)}T_{z_j}^{-1}T_{z_i}\varphi(z),
\end{equation}
where $T_{z_j}^{-1}T_{z_i}$ indicates the shift operator with respect to $z_i\to qz_i$ and $z_j\to q^{-1}z_j$ together. 
(Since we regard $\Phi(z)$ given by (\ref{eq:03Phi0}) as a function of $z=(z_1,z_2,\ldots,z_n) \in(\mathbb{C}^*)^n$ 
under the condition $z_{n+1}=dz_1^{-1}z_2^{-1}\cdots z_n^{-1}$, if we consider the individual shift
$z_i\to qz_i$ ($1\le i\le n$), then the shift $z_{n+1}\to q^{-1}z_{n+1}$ occurs automatically. 
We formally write the shift of this situation by the symbol $T_{z_{n+1}}^{-1}T_{z_i}$ instead of $T_{z_i}$ for our convenience. 
In the same way, we use the symbol $T_{z_i}^{-1}T_{z_{n+1}}$ instead of $T_{z_i}^{-1}$.) 
In this section we use the symbol ${\cal A}$ as the skew-symmetrization with respect to $
S_{n+1}$ on $n+1$ variables $z_1,z_2,\ldots,z_{n+1}$, i.e.,
${\cal A}f(z)=\sum_{\sigma\in S_{n+1}}(\sgn\,\sigma)\,\sigma f(z)$.
Then we have 
\begin{lem}
\label{lem:03nabla=0}
For a function $\varphi(z)$ of $z=(z_1,z_2,\ldots,z_n) \in(\mathbb{C}^*)^n$, 
\begin{equation}
\label{eq:03nabla=0}
\int_0^{\mbox{\small $x$}\8}\Phi(z)\nabla_{\!i,j}\varphi(z)\frac{d_qz_1}{z_1}\wedge\cdots\wedge\frac{d_qz_n}{z_n}=0
\quad (1\le i,j\le n+1), 
\end{equation}
if the integral converges. Moreover, for a  function $\varphi(z)=\varphi(z_1,z_2,\ldots,z_{n+1})$ with 
$z_{n+1}=dz_1^{-1}z_2^{-1}\cdots z_n^{-1}$, 
\begin{equation}
\label{eq:03A1}
\int_0^{\mbox{\small $x$}\8}\Phi(z){\cal A}\nabla_{\!i,j}\varphi(z)\frac{d_qz_1}{z_1}\wedge\cdots\wedge\frac{d_qz_n}{z_n}=0.
\end{equation}
\end{lem}
{\bf Proof.} 
If $i=n+1$ or $j=n+1$, then $T_{z_i}^{-1}T_{z_{n+1}}=T_{z_i}^{-1}$ or $T_{z_{n+1}}^{-1}T_{z_i}=T_{z_i}$, respectively,
so that we can confirm (\ref{eq:03nabla=0}) in the same way as (\ref{eq:00nabla=0}) in Lemma \ref{lem:00nabla=0}. 
While if $i\ne n+1$ and  $j\ne n+1$, 
then from the definition (\ref{eq:03nabla}) of $\nabla_{\!i,j}$, (\ref{eq:03nabla=0}) is equivalent to 
$$
\int_0^{\mbox{\small $x$}\8}\varphi(z)\Phi(z)\frac{d_qz_1}{z_1}\wedge\cdots\wedge\frac{d_qz_n}{z_n}
=\int_0^{\mbox{\small $x$}\8}T_{z_j}^{-1}T_{z_i}\varphi(z)\, T_{z_j}^{-1}T_{z_i}\Phi(z)
\frac{d_qz_1}{z_1}\wedge\cdots\wedge\frac{d_qz_n}{z_n},
$$
which is just confirmed from the fact that the Jackson integral is invariant under the $q$-shift $z_i\to qz_i$ ($1\le i\le n$). 
Next we will confirm (\ref{eq:03A1}). For $\varphi(z)=\varphi(z_1,z_2,\ldots,z_{n+1})$, we have 
$$
\Phi(z){\cal A}\nabla_{\!i,j}\varphi(z)=\Phi(z)\sum_{\sigma\in S_{n+1}}(\sgn\,\sigma)\,\sigma\nabla_{\!i,j}\varphi(z)
=\sum_{\sigma\in S_{n+1}}(\sgn\,\sigma)\,\Phi(z)\nabla_{\!\sigma(i),\sigma(j)}\sigma\varphi(z),
$$
and thus
\begin{eqnarray*}
&&\int_0^{\mbox{\small $x$}\8}\Phi(z){\cal A}\nabla_{\!i,j}\varphi(z)
\frac{d_qz_1}{z_1}\wedge\cdots\wedge\frac{d_qz_n}{z_n}\\
&&=\sum_{\sigma\in S_{n+1}}(\sgn\,\sigma)
\int_0^{\mbox{\small $x$}\8}\Phi(z)
\nabla_{\!\sigma(i),\sigma(j)}\,\sigma\varphi(z)
\frac{d_qz_1}{z_1}\wedge\cdots\wedge\frac{d_qz_n}{z_n}.
\end{eqnarray*}
Since $\sigma\varphi(z)$ is also a function on $(\mathbb{C}^*)^n$ under the condition 
$z_{n+1}=dz_1^{-1}z_2^{-1}\cdots z_n^{-1}$, 
using (\ref{eq:03nabla=0}), we see all terms in the right-hand side of the above equation are equal to $0$. 
$\square$\\

For $c\in \mathbb{C}^*$ and $z=(z_1,z_2.\ldots,z_{n+1})\in (\mathbb{C}^*)^{n+1}$ 
we denote $e(c;z)$ the symmetric polynomial of degree $n+1$ defined by 
$$
e(c;z):=\prod_{i=1}^{n+1}(1-c^{-1}z_i),
$$
(cf.~(\ref{eq:02e(c;z)})) which has a property that $e(c;z)$ vanishes if $z_i=c$ for some $i=1,2,\ldots,n+1$. 
%
%
%
\begin{lem}
\label{lem:03q-diff2}
If we put $\varphi(z)$ as 
\begin{equation}
\label{eq:03varphi0}
\varphi(z)=
\prod_{i=1}^{n+1}(1- a_i^{-1}z_1)(1- b_iz_{n+1})
\times z_2\cdots z_n\prod_{1\le j<k\le n}(z_k-a_j),
\end{equation}
then ${\cal A}\nabla_{\!1,n+1}\varphi(z)
$ is expanded as 
\begin{equation}
\label{eq:03A-nabla-phi0}
{\cal A}\nabla_{\!1,n+1}\varphi(z)=\Big(c_0e_0(z)+c_1e(a_1;z)\Big)\Delta(z),
\end{equation}
where the function $e_0(z)$ is defined by
$$
e_0(z):=1-z_1z_2\cdots z_{n+1}a_1^{-1}a_2^{-1}\cdots a_{n+1}^{-1}
$$
and the constants $c_0$ and $c_1$ are given by 
\begin{equation}
\label{eq:03c0c1}
c_0=(-1)^{n+1}2a_1^nb_1\cdots b_{n+1}\prod_{i=1}^{n+1}(1- a_1^{-1}b_i^{-1}),
\quad c_1=(-1)^{n}2a_1^nb_1\cdots b_{n+1}(1-\prod_{i=1}^{n+1}a_i^{-1}b_i^{-1}).
\end{equation}
In particular, under the condition $z_1z_2\cdots z_{n+1}=d$, the function $e_0(z)$ is a constant and then
\begin{equation}
\label{eq:03c'0}
c'_0:=c_0e_0(z)=(-1)^{n+1} 2a_1^nb_1\cdots b_{n+1}(1-d\,a_1^{-1}a_2^{-1}\cdots a_{n+1}^{-1})\prod_{i=1}^{n+1}(1- a_1^{-1}b_i^{-1}).
\end{equation}
\end{lem}
{\bf Proof.}
Since the ratio $T_{z_{n+1}}^{-1}T_{z_1}\Phi(z)/\Phi(z)$ is written as 
$$
\frac{T_{z_{n+1}}^{-1}T_{z_1}\Phi(z)}{\Phi(z)}
=\prod_{j=1}^{n+1}\frac{(1-b_jz_1)(1-a_j^{-1}z_{n+1})}{(1-qa_j^{-1}z_1)(1-q^{-1}b_jz_{n+1})},
$$
for $\varphi(z)$ given by (\ref{eq:03varphi0}), using (\ref{eq:03nabla}) gives
\begin{equation}
\label{eq:03nabla-phi0}
\nabla_{\!1,n+1}\varphi(z)=
\Big(
\prod_{i=1}^{n+1}(1- a_i^{-1}z_1)(1- b_iz_{n+1})
-\prod_{i=1}^{n+1}(1- b_iz_1)(1- a_i^{-1}z_{n+1})\Big)
\times z_2\cdots z_n\prod_{1\le j<k\le n}(z_k-a_j). 
\end{equation}

Let us use $\nabla$ instead of $\nabla_{\!1,n+1}$ for abbreviation. 
Taking account of the degree of $\nabla\varphi(z)$ as a polynomial of $z$, 
we can expand the skew-symmetrization ${\cal A}\nabla\varphi(z)$ as 
\begin{equation}
\label{eq:03A-nabla-phi1}
{\cal A}\nabla\varphi(z)=\Big(c_0e_0(z)+\sum_{i=1}^nc_ie(a_i;z)+c_{n+1}e(b_{n+1}^{-1};z)\Big)\Delta(z),
\end{equation}
where the coefficients $c_i$ $(i=0,1,\ldots,n)$ are some constants. 
Here we will confirm that $c_i$ vanishes if $i\ge 2$, and $c_0$ and $c_1$ are evaluated as (\ref{eq:03c0c1}).

First we take $z_1=a_1, z_2=a_2,\ldots, z_n=a_n$ and $z_{n+1}=b_{n+1}^{-1}$. 
Then, from  (\ref{eq:03A-nabla-phi1}) with the vanishing property of $e(a_i;z)$, 
we have 
$
{\cal A}\nabla\varphi(z)$ as 
\begin{eqnarray}
\label{eq:03A-nabla-phi1-1}
{\cal A}\nabla\varphi(z)
&=&c_0e_0(z)\Delta(z)
=c_0(1-a_{n+1}^{-1}b_{n+1}^{-1})\prod_{1\le j<k\le n}(a_k-a_j)
\prod_{i=1}^n(b_{n+1}^{-1}-a_i)\nonumber\\
&=&(-1)^{n}c_0a_1\cdots a_n\prod_{i=1}^{n+1}(1-a_i^{-1}b_{n+1}^{-1})\prod_{1\le j<k\le n}(a_k-a_j).
\end{eqnarray}
On the other hand, from (\ref{eq:03nabla-phi0}), we have 
\begin{equation}
\label{eq:03A-nabla-phi1-2}
{\cal A}\nabla\varphi(z)=
-2\prod_{i=1}^{n+1}(1- b_ia_1)(1- a_i^{-1}b_{n+1}^{-1})\times
a_2\cdots a_n\prod_{1\le j<k\le n}(a_k-a_j).
\end{equation}
Comparing (\ref{eq:03A-nabla-phi1-1}) with (\ref{eq:03A-nabla-phi1-2}), 
we obtain $c_0$ in the expression (\ref{eq:03c0c1}).

Next we take $z_1=a_1$, $z_2=a_2,\ldots, z_{n}=a_{n}$ and $z_{n+1}\not\in\{a_1,a_2,\ldots,a_{n}, b_{n+1}^{-1}\}$. Then,
from (\ref{eq:03A-nabla-phi1}) and the vanishing property of $e(a_i;z)$, we have
\begin{equation}
\label{eq:03A-nabla-phi2}
{\cal A}\nabla\varphi(z)=\Big(c_0e_0(z)+c_{n+1} e(b_{n+1}^{-1};z)\Big)\Delta(z).
\end{equation}
On the other hand, from (\ref{eq:03nabla-phi0}) and (\ref{eq:03A-nabla-phi1-2}), we have 
\begin{equation}
\label{eq:03A-nabla-phi3}
{\cal A}\nabla\varphi(z)=
-2\prod_{i=1}^{n+1}(1- b_ia_1)(1- a_i^{-1}z_{n+1})\times
a_2\cdots a_n\prod_{1\le j<k\le n}(a_k-a_j)=c_0e_0(z)\Delta(z).
\end{equation}
Comparing (\ref{eq:03A-nabla-phi2}) and (\ref{eq:03A-nabla-phi3}), we obtain $c_{n+1}=0$. 

Third we take $z_1=a_1$, $z_2=a_2,\ldots, z_{n-1}=a_{n-1}$, $z_{n}\not\in\{a_1,a_2,\ldots,a_{n},b_{n+1}^{-1}\}$ 
and $z_{n+1}=b_{n+1}^{-1}$. Then,
from (\ref{eq:03A-nabla-phi1}) and the vanishing property of $e(a_i;z)$, we have
\begin{equation}
\label{eq:03A-nabla-phi4}
{\cal A}\nabla\varphi(z)=\Big(c_0e_0(z)+c_n e(a_n;z)\Big)\Delta(z).
\end{equation}
On the other hand, from (\ref{eq:03nabla-phi0}), we have 
\begin{equation}
\label{eq:03A-nabla-phi5}
{\cal A}\nabla\varphi(z)=c_0e_0(z)\Delta(z).
\end{equation}
Comparing (\ref{eq:03A-nabla-phi4}) and (\ref{eq:03A-nabla-phi5}), we obtain $c_n=0$. In the same manner, 
by symmetry of ${\cal A}\nabla\varphi(z)$ of $\nabla\varphi(z)$ in (\ref{eq:03nabla-phi0}), 
we also obtain $c_i=0$ if $i\ge 2$. Thus we obtain the expansion (\ref{eq:03A-nabla-phi0}).

Lastly, we take $(z_1,z_2,\ldots,z_{n+1})=(b_1^{-1},b_2^{-1},\ldots,b_{n+1}^{-1})=b^{-1}$. From (\ref{eq:03nabla-phi0}), we have
${\cal A}\nabla\varphi(b^{-1})=0$. On the other hand, from  (\ref{eq:03A-nabla-phi1}), we have 
\begin{equation*}
{\cal A}\nabla\varphi(z)=\Big(c_0 e_0(b^{-1})+c_1 e(a_1;b^{-1})\Big)\Delta(z)
\end{equation*}
so that 
$$c_1=-c_0\frac{e_0(b^{-1})}{e(a_1,b^{-1})}
=(-1)^{n} 2a_1^nb_1\cdots b_{n+1}(1-\prod_{i=1}^{n+1}a_i^{-1}b_i^{-1}),$$
which is what we want to show. 
$\square$\\[10pt]
\noindent
{\bf Proof of Lemma \ref{lem:03q-diff1}.} We will prove (\ref{eq:03rec01}) for $T_{a_j}$ first. 
Without loss of generality, it suffices to show that 
\begin{equation}
\label{eq:03rec0}
T_{a_1}K(x)=K(x)
\frac{1-d\prod_{i=1}^{n+1}a_i^{-1}}
{1-\prod_{i=1}^{n+1}a_i^{-1}b_i^{-1}}
\prod_{i=1}^{n+1}(1-b_i^{-1}a_1^{-1}).
\end{equation}
Since we have 
$
T_{a_1}\Phi(z)/\Phi(z)=
\prod_{i=1}^{n+1}(1-a_1^{-1}z_i)=e(a_1;z)
$
by definition, $T_{a_1}K(z)$ is expressed by 
$$
T_{a_1}K(x)=\int_0^{\mbox{\small $x$}\8}e(a_1;z)\Phi(z)\Delta(z)\,\frac{d_qz_1}{z_1}\wedge\cdots\wedge\frac{d_qz_n}{z_n}.
$$
Applying (\ref{eq:03A1}) in Lemma \ref{lem:03nabla=0} to the fact (\ref{eq:03A-nabla-phi0}) in Lemma \ref{lem:03q-diff2}, we obtain the relation
$$c'_0K(z)+c_1T_{a_1}K(z)=0,$$
where $c_1$ and $c'_0$ are given in (\ref{eq:03c0c1}) and (\ref{eq:03c'0}). This relation coincides with (\ref{eq:03rec0}).

Next we will show the $q$-difference equation (\ref{eq:03rec02}) for $T_{b_j}$ of the case $j=1$ as the same manner as above. 
Since we have 
$
T_{b_1}\Phi(z)/\Phi(z)=
\prod_{i=1}^{n+1}(1-b_1z_i)=e(b_1^{-1};z)$,
$T_{b_1}K(z)$ is expressed by 
$$
T_{b_1}K(x)=\int_0^{\mbox{\small $x$}\8}e(b_1^{-1};z)\Phi(z)\Delta(z)\,\frac{d_qz_1}{z_1}\wedge\cdots\wedge\frac{d_qz_n}{z_n}.
$$
Here if we exchange $a_i$ with $b_i^{-1}$ $(i=1,2,\ldots,n+1)$ 
in the above proof of (\ref{eq:03rec0})
including that of Lemma \ref{lem:03q-diff2}, 
the way of argument is completely symmetric for this exchange. 
Therefore (\ref{eq:03rec02}) for $T_{b_1}$ is obtained exchanging $a_i$ with $b_i^{-1}$ on the coefficient of  (\ref{eq:03rec0}). 

Lastly we will confirm the $q$-difference equations (\ref{eq:03rec03}) and (\ref{eq:03rec04}) for ${\cal K}(x)$. 
From (\ref{eq:03h(x)1}), we have 
$$
T_{a_j}h(z)=h(z)
\quad\mbox{and}\quad 
T_{b_j}h(z)=d(-b_j)^{n+1}h(z).
$$
We therefore obtain (\ref{eq:03rec03}) and (\ref{eq:03rec04})
from (\ref{eq:03rec01}) and (\ref{eq:03rec02}), respectively,
using the above equations. $\square$
%
\subsection{A remark on a relation to the Macdonald identity for $A_n^{(1)}$}
\label{subsection:03-3}
We can use the recurrence relations (\ref{eq:03rec03}) and (\ref{eq:03rec04}) 
directly to the regularization ${\cal K}(z)$ for its evaluation.  
Actually, we can immediately obtain Theorem \ref{thm:03main} as a corollary of Lemma \ref{lem:03q-diff1}.
\begin{cor} 
\label{cor:03cK(x)cK0}
The sum ${\cal K}(x)$ {\rm (}expressed as {\rm (\ref{eq:03cal K(x)1}))} can be written
\begin{equation}
\label{eq:03cal K(x)4}
{\cal K}(x)=
{\cal K}_0(x)\frac{(qa_1^{-1}\cdots a_{n+1}^{-1}d)_\8
(qb_1^{-1}\cdots b_{n+1}^{-1}d^{-1})_\8}
{(qa_1^{-1}\cdots a_{n+1}^{-1}b_1^{-1}\cdots b_{n+1}^{-1})_\8}
\prod_{i=1}^{n+1}\prod_{j=1}^{n+1}(qa_i^{-1}b_j^{-1})_\8,
\end{equation}
where ${\cal K}_0(x)$ is the Jackson integral defined by
\begin{equation}
\label{eq:03cal K0(x)}
{\cal K}_0(x):=\int_0^{\mbox{\small $x$}\8}\frac{1}
{\prod_{1\le i<j\le n+1}(qz_i/z_j)_\8(qz_j/z_i)_\8}\,\frac{d_qz_1}{z_1}\wedge\cdots\wedge\frac{d_qz_n}{z_n}.
\end{equation}
We remark that ${\cal K}_0(x)$ satisfies the 
Macdonald identity for $A_n^{(1)}$ 
in the form proved by Milne {\rm (}straightforward rewrite of {\rm \cite[Theorem 1.58]{Mi85}}{\rm )}, 
i.e., $${\cal K}_0(x)=(1-q)^n(q)_\8^n.$$
\end{cor}
{\bf Proof.} 
Let $T^N$ be the $q$-shift operator with respect to 
$
a_i\to q^{-N}a_i$ and $\ b_i\to q^{-N}b_i\ (i=1,2,\ldots,n+1). 
$
If we set $C$ the right-hand side of (\ref{eq:03cal K(x)4}) without factor ${\cal K}_0(x)$, 
then we have $T^NC\to1$ $(N\to +\8)$. 
From (\ref{eq:03rec03}) and (\ref{eq:03rec04}), 
${\cal K}(x)$ and $C$ satisfy the same recurrence relations, 
so that the ratio ${\cal K}(x)/C$ is invariant under the $q$-shift $T^N$. Therefore we obtain 
$$
\frac{{\cal K}(x)}{C}=T^{N}\frac{{\cal K}(x)}{C}=\lim_{N\to +\8}\frac{T^{N}{\cal K}(x)}{T^{N}C}
=\lim_{N\to +\8}T^N{\cal K}(x),
$$
where 
$$\lim_{N\to +\8}T^N{\cal K}(x)
=
\lim_{N\to +\8}
\int_0^{\mbox{\small $x$}\8}\frac{\prod_{i=1}^{n+1}\prod_{j=1}^{n+1} (q^{1+N}a_j^{-1}z_i)_\8(q^{1+N}b_j^{-1}z_i^{-1})_\8}
{\prod_{1\le i<j\le n+1}(qz_i/z_j)_\8(qz_j/z_i)_\8}\,\frac{d_qz_1}{z_1}\wedge\cdots\wedge\frac{d_qz_n}{z_n},
$$
which coincides with (\ref{eq:03cal K0(x)}) of ${\cal K}_0(x)$. $\square$\\

Though (\ref{eq:03cal K(x)4}) of Corollary \ref{cor:03cK(x)cK0} shows the connection between ${\cal K}(x)$ and ${\cal K}_0(x)$,  
it is not our intention to evaluate ${\cal K}(x)$ via this connection using the Macdonald identity. 
In keeping with our viewpoint our aim is to evaluate ${\cal K}(a)$ directly 
calculating the asymptotic behavior of the truncated Jackson integral. 
\subsection{Asymptotic behavior (Proof of Lemma \ref{lem:03main})}
\label{subsection:03-4}
In this subsection we will give a proof for Lemma \ref{lem:03main}. 
In order to calculate ${\cal K}(a)$, 
we deform $\Phi(z)$ to $\tilde\Phi(z)$ according to
\begin{equation}
\label{eq:03PkP}
\Phi(z)=k(z)\tilde\Phi(z),
\end{equation}
where
\begin{equation}
\label{eq:03k(z)}
k(z):=\prod_{j=1}^{n+1}z_{n+1}^{1-\a_j-\b_j}
\frac{\theta(qa_j^{-1}z_{n+1})}{\theta(b_jz_{n+1})}
=\prod_{j=1}^{n+1}
\frac{\theta(a_jd^{-1}z_1z_2\cdots z_n)}
{(d^{-1}z_1z_2\cdots z_n)^{1-\a_j-\b_j}\theta(qb_j^{-1}d^{-1}z_1z_2\cdots z_n)}
\end{equation}
and
\begin{eqnarray}
\label{eq:03tPhi}
\tilde\Phi(z)&:=&
(z_{n+1}^{-1})^{n+1-\alpha_1-\cdots-\alpha_{n+1}-\beta_1-\cdots-\beta_{n+1}}
\nonumber\\
&&\times\prod_{i=1}^n
\prod_{j=1}^{n+1}\frac{(qa_j^{-1}z_i)_\8}{(b_jz_i)_\8}\times
\prod_{j=1}^{n+1}
\frac{(qb_j^{-1}z_{n+1}^{-1})_\8}{(a_jz_{n+1}^{-1})_\8}
\nonumber\\
&=&
(d^{-1}z_1z_2\cdots z_n)^{n+1-\alpha_1-\cdots-\alpha_{n+1}-\beta_1-\cdots-\beta_{n+1}}
\nonumber\\
&&\times\prod_{i=1}^n
\prod_{j=1}^{n+1}\frac{(qa_j^{-1}z_i)_\8}{(b_jz_i)_\8}\times
\prod_{j=1}^{n+1}
\frac{(qb_j^{-1}d^{-1}z_1\cdots z_n)_\8}{(a_jd^{-1}z_1\cdots z_n)_\8}.
\end{eqnarray}
Here $\alpha_i$ and $\beta_i$ are given by $a_i=q^{\alpha_i}, b_i=q^{\beta_i}$.
Then we define the Jackson integral $\tilde K(x)$ by
\begin{equation}
\label{eq:03tK}
\tilde K(x):=\int_0^{\mbox{\small $x$}\8}
{\tilde\Phi}(z)\Delta(z)\,\frac{d_qz_1}{z_1}\wedge\cdots\wedge\frac{d_qz_n}{z_n}.
\end{equation}
Since $k(z)$ is invariant under the $q$-shift $z_i\to qz_i$, 
from (\ref{eq:03PkP}) and (\ref{eq:03tK}), we have 
\begin{equation}
\label{eq:03KkK}
K(x)=k(x){\tilde K}(x).
\end{equation}
We will calculate the asymptotic behavior of ${\tilde K}(x)$ instead of $K(x)$ when $x=a$. 
Since $\Delta(z)$ is written as 
$$
\Delta(z)=\prod_{i=1}^n(z_{n+1}-z_i)\prod_{1\le j<k\le n}(z_k-z_j)
=\prod_{i=1}^n\frac{1-z_i(d^{-1}z_1\cdots z_{n})}{d^{-1}z_1\cdots z_{n}}\prod_{1\le j<k\le n}(z_k-z_j),
$$
from (\ref{eq:03tPhi}), we have 
\begin{eqnarray*}
\label{eq:03Phi2Delta}
&&\tilde\Phi(z)\Delta(z)=
(d^{-1}z_1z_2\cdots z_n)^{1-\alpha_1-\cdots-\alpha_{n+1}-\beta_1-\cdots-\beta_{n+1}}
\prod_{1\le j<k\le n}(z_k-z_j)\\
&&\times
\prod_{i=1}^n\Bigg(
(1-\frac{z_i}{d}z_1z_2\cdots z_{n})
\prod_{j=1}^{n+1}\frac{(qa_j^{-1}z_i)_\8}{(b_jz_i)_\8}\Bigg)\times
\prod_{j=1}^{n+1}
\frac{(qb_j^{-1}d^{-1}z_1\cdots z_n)_\8}{(a_jd^{-1}z_1\cdots z_n)_\8}.
\end{eqnarray*}
For an integer $N$, let $T^N$ be the $q$-shift operator for a special direction, 
$$
T^N: b_i\to q^{-nN}b_i\ (i=1,2,\ldots,n+1);\ a_j \to q^{(n+1)N}a_j\ (j=1,2,\ldots,n);\ a_{n+1}\to q^{-nN}a_{n+1}. 
$$
Then, by definition $T^N \tilde K(a)$ is written explicitly as
\begin{eqnarray*}
&&T^N \tilde K(a)
=(1-q)^n\!\!\!\!
\sum_{(\nu_1,\ldots,\nu_n)\in \mathbb{N}^n} 
(d^{-1}a_1a_2\cdots a_nq^{\nu_1+\cdots+\nu_n+n(n+1)N})
^{1-\alpha_1-\cdots-\alpha_{n+1}-\beta_1-\cdots-\beta_{n+1}+nN} 
\nonumber\\
&&\times\prod_{i=1}^n\Bigg[
(1-\frac{a_iq^{\nu_i}}{d}a_1a_2\cdots a_{n}q^{\nu_1+\cdots+\nu_n+(n+1)^2N})
\frac{(a_{n+1}^{-1}a_iq^{1+\nu_i+(2n+1)N})_\8}{(b_{n+1}a_iq^{\nu_i+N})_\8}
\prod_{j=1}^{n}\frac{(a_j^{-1}a_iq^{1+\nu_i})_\8}{(b_ja_i^{\nu_i+N})_\8}\Bigg]
\nonumber\\
&&\times
\frac{\prod_{j=1}^{n+1}(b_{j}^{-1}d^{-1}a_1a_2\cdots a_nq^{1+\nu_1+\cdots+\nu_n+n(n+2)N})_\8}
{(a_{n+1}d^{-1}a_1a_2\cdots a_nq^{\nu_1+\cdots+\nu_n+n^2N})_\8
\prod_{j=1}^{n}(a_jd^{-1}z_1z_2\cdots z_nq^{\nu_1+\cdots+\nu_n+(n+1)^2N})_\8}
\\
&&\times
\prod_{1\le i<j\le n}q^{(n+1)N}(a_jq^{\nu_j}-a_iq^{\nu_i}),
\end{eqnarray*}
so that the leading term of the asymptotic behavior of $T^N \tilde K(a)$ as $N\to +\8$ is given by
the term corresponding to $(\nu_1,\ldots,\nu_n)=(0,\ldots,0)$ in the above sum, i.e.,
\begin{eqnarray}
\label{eq:03TNbK(a)}
T^N \tilde K(a)
&\sim& (1-q)^n(d^{-1}a_1a_2\cdots a_nq^{n(n+1)N})
^{1-\alpha_1-\cdots-\alpha_{n+1}-\beta_1-\cdots-\beta_{n+1}+nN} 
\nonumber\\
&&\times\prod_{i=1}^n\prod_{j=1}^{n}(qa_j^{-1}a_i)_\8
\prod_{1\le i<j\le n}q^{(n+1)N}(a_j-a_i)
\nonumber\\
&=& 
(d^{-1}a_1a_2\cdots a_nq^{n(n+1)N})
^{1-\alpha_1-\cdots-\alpha_{n+1}-\beta_1-\cdots-\beta_{n+1}+nN}
\nonumber\\
&&\times (1-q)^n(q)_\8^n \prod_{1\le i<j\le n}q^{(n+1)N} a_j\theta(a_i/a_j)
\quad\quad(N\to +\8).
\end{eqnarray}
On the other hand, if we set $C$ as the right-hand side of (\ref{eq:03cal K(a)}) of Lemma \ref{lem:03main}, 
then we have 
from (\ref{eq:03h(a)2}) and (\ref{eq:03k(z)})
\begin{eqnarray*}
\frac{h(a)C}{k(a)}&=&
(d^{-1}a_1a_2\cdots a_n)^{1-\alpha_1-\cdots-\alpha_{n+1}-\beta_1-\cdots-\beta_{n+1}}
\prod_{1\le i<j\le n}a_j\theta(a_i/a_j)\\
&&
\times\frac{(1-q)^n(q)_\8^n(qb_1^{-1}\cdots b_{n+1}^{-1}d^{-1})_\8\prod_{j=1}^{n+1}(qa_{n+1}^{-1}b_j^{-1})_\8}
{(a_1\cdots a_{n+1}d^{-1})_\8(qa_1^{-1}\cdots a_{n+1}^{-1}b_1^{-1}\cdots b_{n+1}^{-1})_\8\prod_{i=1}^{n}\prod_{j=1}^{n+1}(b_ja_i)_\8},
\end{eqnarray*}
and thus  
\begin{eqnarray}
\label{eq:03TNhCk}
&&T^N\frac{h(a)C}{k(a)}=
(d^{-1}a_1a_2\cdots a_nq^{n(n+1)N})^{1-\alpha_1-\cdots-\alpha_{n+1}-\beta_1-\cdots-\beta_{n+1}+nN}
\prod_{1\le i<j\le n}q^{(n+1)N}a_j\theta(a_i/a_j)
\nonumber\\
&&
\quad\times\frac{(1-q)^n(q)_\8^n
(b_1^{-1}\cdots b_{n+1}^{-1}d^{-1}q^{1+n(n+1)N})_\8\prod_{j=1}^{n+1}(a_{n+1}^{-1}b_j^{-1}q^{1+2nN})_\8}
{(a_1\cdots a_{n+1}d^{-1}q^{n^2N})_\8(a_1^{-1}\cdots a_{n+1}^{-1}b_1^{-1}\cdots b_{n+1}^{-1}q^{1+nN})_\8\prod_{i=1}^{n}\prod_{j=1}^{n+1}(b_ja_iq^N)_\8}.
\end{eqnarray}
Since ${\cal K}(a)$ and $C$ satisfy the same recurrence relations with respect to $a_i$ and $b_i$, 
the ratio ${\cal K}(a)/C$ is invariant under the $q$-shift $T^N$. From (\ref{eq:03TNbK(a)}) and (\ref{eq:03TNhCk}), 
we therefore obtain 
\begin{eqnarray*}
&&\frac{{\cal K}(a)}{C}=T^N\frac{{\cal K}(a)}{C}=\frac{T^N\tilde K(a)}{T^Nh(a)C/k(a)}=
\lim_{N\to +\8}\frac{T^N\tilde K(a)}{T^Nh(a)C/k(a)}\\
&&=\lim_{N\to +\8}
\frac{(a_1\cdots a_{n+1}d^{-1}q^{n^2N})_\8
(a_1^{-1}\cdots a_{n+1}^{-1}b_1^{-1}\cdots b_{n+1}^{-1}q^{1+nN})_\8
\prod_{i=1}^{n}\prod_{j=1}^{n+1}(b_ja_iq^N)_\8}
{(b_1^{-1}\cdots b_{n+1}^{-1}d^{-1}q^{1+n(n+1)N})_\8\prod_{j=1}^{n+1}(a_{n+1}^{-1}b_j^{-1}q^{1+2nN})_\8}\\
&&=1,
\end{eqnarray*}
which completes the proof of Lemma \ref{lem:03main}. $\square$

\subsection*{Acknowledgements}
This work was supported by the Australian Research Council and JSPS KAKENHI Grant Number 25400118.

\end{document}